\documentclass[11pt]{article}

\usepackage[latin1]{inputenc}\usepackage{epsfig}\usepackage{amsfonts}

\title{
  {\huge Fourier Theory on the Complex Plane I} \\
  Conjugate Pairs of Fourier Series \\
  and Inner Analytic Functions }

\author{
  \Large Jorge L. deLyra \\
  Department of Mathematical Physics \\
  Physics Institute \\
  University of São Paulo }

\date{March 2, 2015}

\newcommand{\ii}{\mbox{\boldmath$\imath$}}

\newcommand{\FFrac}[2]{{\displaystyle\frac{\displaystyle #1}{\displaystyle #2}}}

\newcommand{\at}[2]{\left.\rule{0em}{3ex}\right[_{\,#1}^{\,#2}}

\newcommand{\e}[1]{\,{\rm e}^{#1}}

\begin{document}\maketitle

\begin{abstract}
  \noindent
  A correspondence between arbitrary Fourier series and certain analytic
  functions on the unit disk of the complex plane is established. The
  expression of the Fourier coefficients is derived from the structure of
  complex analysis. The orthogonality and completeness relations of the
  Fourier basis are derived in the same way. It is shown that the limiting
  function of any Fourier series is also the limit to the unit circle of
  an analytic function in the open unit disk. An alternative way to
  recover the original real functions from the Fourier coefficients, which
  works even when the Fourier series are divergent, is thus presented. The
  convergence issues are discussed up to a certain point. Other possible
  uses of the correspondence established are pointed out.
\end{abstract}

\section{Introduction}

In this paper we will establish an interesting relation between Fourier
series and analytic functions. This leads to an alternative way to deal
with Fourier series and to characterize the corresponding real functions.
This relation will allow us to discuss the convergence of Fourier series
in terms of the convergence of the Taylor series of analytic functions.
The convergence issues will be developed up to a certain point, and
further developments will be discussed in a follow-up
paper~\cite{FTotCPII}. This relation will also give us an alternative way
to recover, from the coefficients of the series, the functions that
originated them, which works even if the Fourier series are divergent.
Perhaps most importantly, it will provide a different and possibly richer
point of view for Fourier series and the corresponding real functions.

We will use repeatedly the following very well-known and fundamental
theorem of complex analysis, about complex power series~\cite{ChurComp}.
If we consider the general complex power series written around the origin
$z=0$,

\begin{displaymath}
  S_{z}
  =
  \sum_{k=0}^{\infty}
  c_{k}z^{k},
\end{displaymath}

\noindent
where $z=x+\ii y$ is a complex variable and $c_{k}$ are arbitrary complex
constants, then the following holds. If $S_{z}$ converges at a point
$z_{1}\neq 0$, then it is convergent and absolutely convergent on an open
disk centered at $z=0$ with its boundary passing through $z_{1}$. In
addition to this, it converges uniformly on any closed set contained
within this open disk. We will refer to this state of affairs in what
regards convergence as {\em strong convergence}, and will refer to this
theorem as the {\em basic convergence theorem}. Furthermore, the power
series converges to an analytic function of which it is the Taylor series
around $z=0$.

We will make a conceptual distinction between trigonometric series and
Fourier series. An arbitrary real trigonometric series on the real
variable $\theta$ with domain in the periodic interval $[-\pi,\pi]$ is
given by

\begin{displaymath}
  S
  =
  \frac{1}{2}\,
  \alpha_{0}
  +
  \sum_{k=1}^{\infty}
  \alpha_{k}\cos(k\theta)
  +
  \sum_{k=1}^{\infty}
  \beta_{k}\sin(k\theta),
\end{displaymath}

\noindent
where $\alpha_{k}$ and $\beta_{k}$ are any real numbers. If there is a
real function $f(\theta)$ such that the coefficients $\alpha_{k}$ and
$\beta_{k}$ are given in terms of that function by the integrals

\noindent
\begin{eqnarray*}
  \alpha_{k}
  & = &
  \frac{1}{\pi}
  \int_{-\pi}^{\pi}d\theta\,
  f(\theta)
  \cos(k\theta),
  \\
  \beta_{k}
  & = &
  \frac{1}{\pi}
  \int_{-\pi}^{\pi}d\theta\,
  f(\theta)
  \sin(k\theta),
\end{eqnarray*}

\noindent
then we call this series the Fourier series of that
function~\cite{ChurFour}. Since the Fourier coefficients are defined by
means of integrals, it is clear that one can add to $f(\theta)$ any
zero-measure function without modifying them. This means that a convergent
Fourier series can only be said to converge {\em almost everywhere} to the
function which originated it, that is, with the possible exclusion of a
zero-measure subset of the domain.

With this limitation, the Fourier basis in the space of real functions on
the periodic interval, formed by the constant function, the set of
functions $\cos(k\theta)$, with $k=1,\ldots,\infty$, and the set of
functions $\sin(k\theta)$, with $k=1,\ldots,\infty$, is complete to
generate all sufficiently well-behaved functions in that interval. With
the exclusion of the constant function, the remaining basis generates the
set of all sufficiently well-behaved zero-average real functions. The
remaining basis functions satisfy the orthogonality relations

\noindent
\begin{eqnarray}\label{EQorthrels}
  \frac{1}{\pi}
  \int_{-\pi}^{\pi}d\theta\,
  \cos(k_{1}\theta)
  \cos(k_{2}\theta)
  & = &
  \delta_{k_{1}k_{2}},
  \nonumber\\
  \frac{1}{\pi}
  \int_{-\pi}^{\pi}d\theta\,
  \sin(k_{1}\theta)
  \sin(k_{2}\theta)
  & = &
  \delta_{k_{1}k_{2}},
  \nonumber\\
  \frac{1}{\pi}
  \int_{-\pi}^{\pi}d\theta\,
  \cos(k_{1}\theta)
  \sin(k_{2}\theta)
  & = &
  0,
\end{eqnarray}

\noindent
for $k_{1}\geq 1$ and $k_{2}\geq 1$. In a stricter sense, the
good-behavior conditions on the real functions are that they be integrable
and that they be such that their Fourier series converge. However, one
might consider the set of coefficients themselves to be sufficient to
characterize the function that originated them, even if the series does
not converge. This makes full sense if there is an alternative way to
recover the functions from the coefficients of their series, such as the
one we will present in this paper, which is not dependent on the
convergence of the series. In this case the condition of integrability
suffices. The important point to be kept in mind here is that the set of
coefficients determines the function uniquely almost everywhere over the
periodic interval.

The parity of the real functions will play an important role in this
paper. Any real function defined in a symmetric domain around zero,
without any additional hypotheses, can be separated into its even and odd
parts. An even function is one that satisfies the condition
$f(-\theta)=f(\theta)$, while an odd one satisfies the condition
$f(-\theta)=-f(\theta)$. For any real function we can write that
$f(\theta)=f_{\rm c}(\theta)+f_{\rm s}(\theta)$, with

\noindent
\begin{eqnarray*}
  f_{\rm c}(\theta)
  & = &
  \frac{f(\theta)+f(-\theta)}{2},
  \\
  f_{\rm s}(\theta)
  & = &
  \frac{f(\theta)-f(-\theta)}{2},
\end{eqnarray*}

\noindent
where $f_{\rm c}(\theta)$ is even and $f_{\rm s}(\theta)$ is odd. The
Fourier basis can also be separated into even and odd parts. Since the
constant function and the cosines are even, they generate the even parts
of the real functions, while the sines, being odd, generate the odd parts.
Since the set of cosines and the set of sines are two independent and
mutually orthogonal sets of functions, the convergence of the
trigonometric series can only be accomplished by the separate convergence
of the cosine sub-series and the sine sub-series, which we denote by

\noindent
\begin{eqnarray*}
  S_{\rm c}
  & = &
  \sum_{k=1}^{\infty}
  \alpha_{k}\cos(k\theta),
  \\
  S_{\rm s}
  & = &
  \sum_{k=1}^{\infty}
  \beta_{k}\sin(k\theta),
\end{eqnarray*}

\noindent
so that $S=\alpha_{0}/2+S_{\rm c}+S_{\rm s}$. For simplicity, in this
paper we will consider only series in which $\alpha_{0}=0$, which
correspond to functions $f(\theta)$ that have zero average over the
periodic interval. There is of course no loss of generality involved in
doing this, since the addition of a constant term is a trivial procedure
that does not bear on the convergence issues. The discussion of the
convergence of arbitrary trigonometric series $S$ is therefore equivalent
to the separate discussion of the convergence of arbitrary cosine series
$S_{\rm c}$ and arbitrary sine series $S_{\rm s}$. These last two classes
consist of trigonometric series with definite parities and we will name
them {\em Definite-Parity} trigonometric series, or DP trigonometric
series for short.

A note about the fact that we limit ourselves to real Fourier series here.
We might as well consider the series $S$ with complex coefficients
$\alpha_{k}$ and $\beta_{k}$, but due to the linearity of the series with
respect to these coefficients, any such complex Fourier series would at
once decouple into two real Fourier series, one in the real part and one
in the imaginary part, and would therefore reduce the discussion to the
one we choose to develop here. Therefore nothing fundamentally new is
introduced by the examination of complex Fourier series, and it is enough
to limit the discussion to the real case.

\section{Trigonometric Series on the Complex Plane}

First of all, let us establish a very basic correspondence between real
trigonometric series and power series in the complex plane. In this
section we do not assume that the trigonometric series are Fourier series.
In fact, for the time being we impose no additional restrictions on the
numbers $\alpha_{k}$ and $\beta_{k}$, other than that they be real, and in
particular we do not assume anything about the convergence of the series.

Consider then an arbitrary DP trigonometric series. We now introduce a
useful definition. Given a cosine series $S_{\rm c}$ with coefficients
$\alpha_{k}$, we will define from it a corresponding sine series by

\begin{displaymath}
  \bar{S}_{\rm c}
  =
  \sum_{k=1}^{\infty}
  \alpha_{k}\sin(k\theta).
\end{displaymath}

\noindent
We will call this new trigonometric series the {\em Fourier-Conjugate}
series to $S_{\rm c}$, or the FC series for short. Note that $\bar{S}_{\rm
  c}$ is odd instead of even. Similarly, given a sine series $S_{\rm s}$
with coefficients $\beta_{k}$, we will define from it a corresponding
cosine series by

\begin{displaymath}
  \bar{S}_{\rm s}
  =
  \sum_{k=1}^{\infty}
  \beta_{k}\cos(k\theta),
\end{displaymath}

\noindent
which we will also name the Fourier-Conjugate series to $S_{\rm s}$. Note
that $\bar{S}_{\rm s}$ is even instead of odd. We see therefore that the
set of all DP trigonometric series can be organized in pairs of mutually
conjugate series. In any given pair, each series is the FC series of the
other.

From now on we will denote all trigonometric series coefficients by
$a_{k}$, regardless of whether the series originally under discussion is a
cosine series or a sine series. Now, given any cosine series $S_{\rm c}$
or any sine series $S_{\rm s}$, we may define from it a complex series
$S_{v}$ by the use of the original series and its FC series as the real
and imaginary parts of the complex series. In the case of an original
cosine series we thus define

\begin{displaymath}
  S_{v}
  =
  S_{\rm c}
  +
  \ii
  \bar{S}_{\rm c},
\end{displaymath}

\noindent
while in the case of an original sine series we define

\begin{displaymath}
  S_{v}
  =
  \bar{S}_{\rm s}
  +
  \ii
  S_{\rm s}.
\end{displaymath}

\noindent
In this way the discussion of the convergence of arbitrary DP
trigonometric series can be reduced to the discussion of the convergence
of the corresponding complex series $S_{v}$. In either one of the two
cases above this series may be written as

\noindent
\begin{eqnarray*}
  S_{v}
  & = &
  \sum_{k=1}^{\infty}
  a_{k}\cos(k\theta)
  +
  \ii
  \sum_{k=1}^{\infty}
  a_{k}\sin(k\theta)
  \\
  & = &
  \sum_{k=1}^{\infty}
  a_{k}
  \left[
    \cos(k\theta)
    +
    \ii
    \sin(k\theta)
  \right],
\end{eqnarray*}

\noindent
where the coefficients $a_{k}$ are still completely arbitrary. If we now
define the complex variable $v=\exp(\ii\theta)$, then using Euler's
formula we may write this complex series as

\begin{displaymath}
  S_{v}
  =
  \sum_{k=1}^{\infty}
  a_{k}v^{k},
\end{displaymath}

\noindent
so that it becomes, therefore, a complex power series with real
coefficients on the unit circle centered at the origin, in the complex
plane. Finally, we may look at this series as a restriction to the unit
circle of a full power series on the complex plane if we introduce an
extra real variable $\rho\geq 0$, so that a complex variable

\noindent
\begin{eqnarray*}
  z
  & = &
  \rho v
  \\
  & = &
  \rho\e{\ii\theta}
\end{eqnarray*}

\noindent
can be defined over the whole complex plane, and consider the complex
power series, still with real coefficients,

\begin{displaymath}
  S_{z}
  =
  \sum_{k=1}^{\infty}
  a_{k}z^{k}.
\end{displaymath}

\noindent
This is a complex power series centered at $z=0$, with no $k=0$ term, so
that is assumes the value zero at $z=0$. Apart from the fact that
$a_{0}=0$, it has real but otherwise arbitrary coefficients. The series
$S_{v}$ that we constructed from a pair of FC trigonometric series is just
$S_{z}$ restricted to $\rho=1$, for $\theta\in[-\pi,\pi]$. In other words,
the series $S_{v}$ is a restriction to the unit circle of the complex
power series we just defined. There is, therefore, a one-to-one
correspondence between pairs of mutually FC trigonometric series and
complex power series around $z=0$ with real coefficients and $a_{0}=0$.

We thus establish that the discussion of the convergence of arbitrary DP
trigonometric series can be reduced to the discussion of the convergence
of the corresponding complex power series $S_{z}$ on the unit circle. In
fact, the whole question of the convergence of trigonometric series is
revealed to be identical to the question of the convergence of complex
power series on the boundary of the unit disk, including the cases in
which that disk is the maximum disk of convergence of the power series.

\section{Fourier Series on the Complex Plane}

Let us now show that the usual formulas giving the Fourier coefficients,
in terms of integrals involving the corresponding real functions, follow
as consequences of the analytic properties of certain complex functions
within the open unit disk. In order to do this, let us consider a pair of
FC trigonometric series that have the rather weak property that there is
at least one value of $\theta$ for which both elements of the pair are
convergent. Note that this constitutes an indirect restriction on the
coefficients of the series. It follows at once that the power series
$S_{z}$ converges at the point on the unit circle that corresponds to that
value of $\theta$. Consequently, it follows from the basic convergence
theorem that the power series is strongly convergent at least on the open
unit disk. Furthermore, it converges to a complex function that is
analytic at least on the open unit disk, which we will denote by $w(z)$,
and therefore we may now write

\begin{displaymath}
  w(z)
  =
  \sum_{k=1}^{\infty}
  a_{k}z^{k}.
\end{displaymath}

\noindent
Note that, since the coefficients are real, the function $w(z)$ reduces to
a purely real function on the open interval $(-1,1)$ of the real axis. It
is therefore the analytic continuation of a real analytic function defined
on that interval. Apart from this fact, from the fact that $w(0)=0$, and
from the fact that it is analytic on the open unit disk, it is an
otherwise arbitrary analytic function. In addition to all this, we have
that $S_{z}$ is the Taylor series of $w(z)$ around $z=0$. We will call an
analytic function that has these properties an {\em inner analytic
  function}. Let us list the defining properties. An inner analytic
function is one that:

\begin{itemize}

\item is analytic at least on the open unit disk;

\item is the analytic continuation of a real function defined in $(-1,1)$;

\item assumes the value zero at $z=0$.

\end{itemize}

\noindent
Let us now examine another property of $w(z)$ implied by the fact that it
is the analytic continuation of a real function. If we use polar
coordinates and write $z=\rho\exp(\ii\theta)$, with $\theta\in[-\pi,\pi]$,
then we may write out the Taylor series of $w(z)$ as

\begin{displaymath}
  w(z)
  =
  \sum_{k=1}^{\infty}
  a_{k}\rho^{k}
  \left[
    \cos(k\theta)
    +
    \ii
    \sin(k\theta)
  \right].
\end{displaymath}

\noindent
Since the coefficients are real, we have at once that

\begin{displaymath}
  w(z)
  =
  \left[
    \,
    \sum_{k=1}^{\infty}
    a_{k}\rho^{k}
    \cos(k\theta)
  \right]
  +
  \ii
  \left[
    \,
    \sum_{k=1}^{\infty}
    a_{k}\rho^{k}
    \sin(k\theta)
  \right],
\end{displaymath}

\noindent
where the expressions within square brackets are real. If we write $w(z)$
in terms of its real and imaginary parts,

\begin{displaymath}
  w(z)
  =
  f_{\rm c}(\rho,\theta)
  +
  \ii
  f_{\rm s}(\rho,\theta),
\end{displaymath}

\noindent
then the real part $f_{\rm c}(\rho,\theta)$ must be even on $\theta$,
because it is the function that the cosine series contained in $S_{z}$
converges to,

\begin{displaymath}
  f_{\rm c}(\rho,\theta)
  =
  \sum_{k=1}^{\infty}
  a_{k}\rho^{k}
  \cos(k\theta).
\end{displaymath}

\noindent
Similarly, the imaginary part $f_{\rm s}(\rho,\theta)$ must be odd on
$\theta$, because it is the function that the sine series contained in
$S_{z}$ converges to,

\begin{displaymath}
  f_{\rm s}(\rho,\theta)
  =
  \sum_{k=1}^{\infty}
  a_{k}\rho^{k}
  \sin(k\theta).
\end{displaymath}

\noindent
With these preliminaries established, we may now proceed towards our
objective here, which consists of the inversion of the relations above, so
that we may write $a_{k}$ in terms of $f_{\rm c}(\rho,\theta)$, or in
terms of $f_{\rm s}(\rho,\theta)$, by means of the use of the analytic
structure within the open unit disk. Consider then the Cauchy integral
formulas for the function $w(z)$ and its derivatives, written around $z=0$
for the $k^{\rm th}$ derivative,

\begin{displaymath}
  w^{k\prime}(0)
  =
  \frac{k!}{2\pi\ii}
  \oint_{C}dz\,
  \frac{w(z)}{z^{k+1}},
\end{displaymath}

\noindent
where $C$ is the circle centered at $z=0$ with radius $\rho\in(0,1)$. The
coefficients of the Taylor series of $w(z)$ may be written in terms of
these integrals, so that we have for $a_{k}$

\noindent
\begin{eqnarray*}
  a_{k}
  & = &
  \frac{w^{k\prime}(0)}{k!}
  \\
  & = &
  \frac{1}{2\pi\ii}
  \oint_{C}dz\,
  \frac{w(z)}{z^{k+1}}.
\end{eqnarray*}

\noindent
It is very important to note that since $w(z)$ is analytic in the open
unit disk, by the Cauchy-Goursat theorem the integral is independent of
$\rho$ within that disk, and therefore so are the coefficients $a_{k}$.
We now write the integral explicitly, using the integration variable
$\theta$ on the circle of radius $\rho$. We have
$dz=\ii\rho\exp(\ii\theta)d\theta$, and therefore get

\noindent
\begin{eqnarray*}
  a_{k}
  & = &
  \frac{1}{2\pi\ii}
  \int_{-\pi}^{\pi}d\theta\,
  \ii\rho^{-k}\e{-\ii k\theta}
  w(z)
  \\
  & = &
  \frac{\rho^{-k}}{2\pi}
  \int_{-\pi}^{\pi}d\theta\,
  \left[
    f_{\rm c}(\rho,\theta)
    +
    \ii
    f_{\rm s}(\rho,\theta)
  \right]
  \left[
    \cos(k\theta)
    -
    \ii
    \sin(k\theta)
  \right]
  \\
  & = &
  \frac{\rho^{-k}}{2\pi}
  \int_{-\pi}^{\pi}d\theta\,
  \left\{
    \rule{0em}{2.5ex}
    \left[
      f_{\rm c}(\rho,\theta)
      \cos(k\theta)
      +
      f_{\rm s}(\rho,\theta)
      \sin(k\theta)
    \right]
  \right.
  +
  \\
  &   &
  \hspace{6em}
  +
  \left.
    \ii
    \left[
      f_{\rm s}(\rho,\theta)
      \cos(k\theta)
      -
      f_{\rm c}(\rho,\theta)
      \sin(k\theta)
    \right]
    \rule{0em}{2.5ex}
  \right\}.
\end{eqnarray*}

\noindent
Since we know that $a_{k}$ are real, we may at once conclude that the
imaginary part of this last integral is zero. But we can state more than
just that, because all the functions appearing in all these integrals have
definite parities on $\theta$, and hence we see that the integrands that
appear in the imaginary part are odd, while the integrals are over
symmetric intervals. We therefore conclude that the following two
integrals are separately zero,

\noindent
\begin{eqnarray*}
  \int_{-\pi}^{\pi}d\theta\,
  f_{\rm s}(\rho,\theta)
  \cos(k\theta)
  & = &
  0,
  \\
  \int_{-\pi}^{\pi}d\theta\,
  f_{\rm c}(\rho,\theta)
  \sin(k\theta)
  & = &
  0,
\end{eqnarray*}

\noindent
for all $k$. We are therefore left with the following expression for
$a_{k}$,

\begin{equation}\label{EQcoefsak}
  a_{k}
  =
  \frac{\rho^{-k}}{2\pi}
  \int_{-\pi}^{\pi}d\theta\,
  \left[
    f_{\rm c}(\rho,\theta)
    \cos(k\theta)
    +
    f_{\rm s}(\rho,\theta)
    \sin(k\theta)
  \right].
\end{equation}

\noindent
In order to continue the analysis of the coefficients $a_{k}$ we consider
now the following integral on the same circuit $C$,

\begin{displaymath}
  \oint_{C}dz\,
  w(z)z^{k-1}
  =
  0,
\end{displaymath}

\noindent
with $k>0$. The integral is zero by the Cauchy-Goursat theorem, since for
$k>0$ the integrand is analytic on the open unit disk. As before we write
the integral on the circle of radius $\rho$ using the integration variable
$\theta$, to get

\noindent
\begin{eqnarray*}
  0
  & = &
  \int_{-\pi}^{\pi}d\theta\,
  \ii\rho^{k}\e{\ii k\theta}
  w(z)
  \\
  & = &
  \ii\rho^{k}
  \int_{-\pi}^{\pi}d\theta\,
  \left[
    f_{\rm c}(\rho,\theta)
    +
    \ii
    f_{\rm s}(\rho,\theta)
  \right]
  \left[
    \cos(k\theta)
    +
    \ii
    \sin(k\theta)
  \right]
  \\
  & = &
  \ii\rho^{k}
  \int_{-\pi}^{\pi}d\theta\,
  \left\{
    \rule{0em}{2.5ex}
    \left[
      f_{\rm c}(\rho,\theta)
      \cos(k\theta)
      -
      f_{\rm s}(\rho,\theta)
      \sin(k\theta)
    \right]
  \right.
  +
  \\
  &   &
  \hspace{6em}
  +
  \left.
    \ii
    \left[
      f_{\rm s}(\rho,\theta)
      \cos(k\theta)
      +
      f_{\rm c}(\rho,\theta)
      \sin(k\theta)
    \right]
    \rule{0em}{2.5ex}
  \right\}.
\end{eqnarray*}

\noindent
Once again the integrals that appear in the imaginary part of this last
expression are zero by parity arguments, and since $\rho\neq 0$ we are
left with

\begin{displaymath}
  \int_{-\pi}^{\pi}d\theta\,
  \left[
    f_{\rm c}(\rho,\theta)
    \cos(k\theta)
    -
    f_{\rm s}(\rho,\theta)
    \sin(k\theta)
  \right]
  =
  0,
\end{displaymath}

\noindent
which is valid for all $k>0$. We conclude therefore that the two integrals
shown are equal,

\begin{displaymath}
  \int_{-\pi}^{\pi}d\theta\,
  f_{\rm c}(\rho,\theta)
  \cos(k\theta)
  =
  \int_{-\pi}^{\pi}d\theta\,
  f_{\rm s}(\rho,\theta)
  \sin(k\theta),
\end{displaymath}

\noindent
for all $k>0$. If we now go back to the expression in
Equation~(\ref{EQcoefsak}) for $a_{k}$ we see that the two integrals
appearing in that expression are equal to each other. We may therefore
write for the coefficients

\noindent
\begin{eqnarray*}
  a_{k}
  & = &
  \frac{\rho^{-k}}{\pi}
  \int_{-\pi}^{\pi}d\theta\,
  f_{\rm c}(\rho,\theta)
  \cos(k\theta)
  \\
  & = &
  \frac{\rho^{-k}}{\pi}
  \int_{-\pi}^{\pi}d\theta\,
  f_{\rm s}(\rho,\theta)
  \sin(k\theta).
\end{eqnarray*}

\noindent
We observe now that these formulas for the coefficients $a_{k}$ are simple
extensions of the usual formulas for the Fourier coefficients of the even
function $f_{\rm c}(\rho,\theta)$ and the odd function $f_{\rm
  s}(\rho,\theta)$, and therefore are related in a simple way to the
Fourier coefficients for the real function of $\theta$

\begin{displaymath}
  f(\rho,\theta)
  =
  f_{\rm c}(\rho,\theta)
  +
  f_{\rm s}(\rho,\theta),
\end{displaymath}

\noindent
with $\rho$ interpreted as an extra parameter. In fact, these formulas
become the usual ones in the $\rho\to 1$ limit, thus completing the
construction of a pair of FC Fourier series on the unit circle.

Whether or not we may now take the limit $\rho\to 1$ in these formulas
depends on whether or not the coefficients, and hence the integrals that
define them, are continuous functions of $\rho$ at the unit circle, for
limits coming from {\em within} the unit disk. We saw before that the
coefficients are constant with $\rho$, and therefore are continuous
functions of $\rho$ within the open unit disk. We therefore know that
their $\rho\to 1$ limits exist. Furthermore, by construction these are the
coefficients of the FC pair of DP trigonometric series we started with, on
the unit circle. Therefore the coefficients assume at $\rho=1$ the values
given by their limits when $\rho\to 1$.

Consequently, the coefficients $a_{k}$ and the expressions giving them
within the open unit disk are continuous from within at the unit circle,
as functions of $\rho$, and so are the integrals appearing in those
expressions. We may now take the $\rho\to 1$ limit and therefore get the
usual formulas for the Fourier coefficients,

\noindent
\begin{eqnarray*}
  a_{k}
  & = &
  \alpha_{k}
  \\
  & = &
  \beta_{k},
  \\
  \alpha_{k}
  & = &
  \frac{1}{\pi}
  \int_{-\pi}^{\pi}d\theta\,
  f_{\rm c}(1,\theta)
  \cos(k\theta),
  \\
  \beta_{k}
  & = &
  \frac{1}{\pi}
  \int_{-\pi}^{\pi}d\theta\,
  f_{\rm s}(1,\theta)
  \sin(k\theta),
\end{eqnarray*}

\noindent
where

\noindent
\begin{eqnarray*}
  f_{\rm c}(1,\theta)
  & = &
  \lim_{\rho\to 1}
  f_{\rm c}(\rho,\theta),
  \\
  f_{\rm s}(1,\theta)
  & = &
  \lim_{\rho\to 1}
  f_{\rm s}(\rho,\theta).
\end{eqnarray*}

\noindent
We see therefore that the two trigonometric series of the pair of FC
series we started with, under the very weak hypothesis that they both
converge at one common point, are in fact the DP Fourier series of the DP
functions $f_{\rm c}(1,\theta)$ and $f_{\rm s}(1,\theta)$ which are
obtained as the $\rho\to 1$ limits of the real part $f_{\rm
  c}(\rho,\theta)$ and of the imaginary part $f_{\rm s}(\rho,\theta)$ of
the inner analytic function $w(z)$.

It is important to note that $w(z)$ might not be analytic at some points
on the unit circle. Also, so far we cannot state that the Taylor series
$S_{z}$ converges anywhere on the unit circle, besides that single point
at which we assumed the convergence of the pair of FC trigonometric
series. For it to be possible to define the real integrals over the unit
circle, the $\rho\to 1$ limits of the functions $f_{\rm c}(\rho,\theta)$
and $f_{\rm s}(\rho,\theta)$ must exist at least almost everywhere on the
unit circle parametrized by $\theta$. They may fail to exist at points
where $w(z)$ has isolated singularities on that circle. Therefore, for the
moment the definition of the trigonometric series as Fourier series on the
unit circle must remain conditioned to the existence of these limits
almost everywhere.

Note that, if the limits to the unit circle result in isolated
singularities in $f_{\rm c}(1,\theta)$ or $f_{\rm s}(1,\theta)$, then
these must be integrable ones along the unit circle, since the $a_{k}$
coefficients are all finite.

\section{Fourier-Taylor Correspondence}

Let us now show that there is a complete one-to-one correspondence between
arbitrarily given pairs of real FC Fourier series and the inner analytic
functions within the open unit disk. To this end, let us imagine that one
begins the whole argument of the last section over again, but this time
starting with a pair of FC Fourier series. What this means is that there
is a zero-average real function $f(\theta)$ defined in the periodic
interval such that the coefficients $a_{k}$ are given in terms of that
function by the integrals

\noindent
\begin{eqnarray}\label{EQfrcoefs}
  a_{k}
  & = &
  \frac{1}{\pi}
  \int_{-\pi}^{\pi}d\theta\,
  f_{\rm c}(\theta)
  \cos(k\theta)
  \nonumber\\
  & = &
  \frac{1}{\pi}
  \int_{-\pi}^{\pi}d\theta\,
  f_{\rm s}(\theta)
  \sin(k\theta),
\end{eqnarray}

\noindent
where we have $f(\theta)=f_{\rm c}(\theta)+f_{\rm s}(\theta)$, with
$f_{\rm c}(\theta)$ even on $\theta$ and $f_{\rm s}(\theta)$ odd on
$\theta$. The two Fourier series generated by $f_{\rm c}(\theta)$ and
$f_{\rm s}(\theta)$ have exactly the same coefficients, and are therefore
the FC series of one another. We may therefore consider the corresponding
functions to be FC functions of one another as well, even if the series do
not converge. We will denote these FC functions by $\bar{f}_{\rm
  c}(\theta)$ in the case of an original cosine series with coefficients
given by $f_{\rm c}(\theta)$, and by $\bar{f}_{\rm s}(\theta)$ in the case
of an original sine series with coefficients given by $f_{\rm
  s}(\theta)$. We have therefore that

\noindent
\begin{eqnarray*}
  f_{\rm c}(\theta)
  & = &
  \bar{f}_{\rm s}(\theta),
  \\
  f_{\rm s}(\theta)
  & = &
  \bar{f}_{\rm c}(\theta).
\end{eqnarray*}

\noindent
Note that $\bar{f}_{\rm c}(\theta)$ is in fact odd instead of even, while
$\bar{f}_{\rm s}(\theta)$ is in fact even instead of odd.

As before, we assume that there is at least one value of $\theta$ for
which both series in the FC pair converge. We may then use the
coefficients $a_{k}$ to define the inner analytic function $w(z)$, as we
did in the previous section, we may identify these coefficients as those
of its Taylor series, and therefore these same coefficients turn out to be
given by

\noindent
\begin{eqnarray*}
  a_{k}
  & = &
  \frac{\rho^{-k}}{\pi}
  \int_{-\pi}^{\pi}d\theta\,
  f_{\rm c}(\rho,\theta)
  \cos(k\theta)
  \\
  & = &
  \frac{\rho^{-k}}{\pi}
  \int_{-\pi}^{\pi}d\theta\,
  f_{\rm s}(\rho,\theta)
  \sin(k\theta),
\end{eqnarray*}

\noindent
where $f_{\rm c}(\rho,\theta)$ and $f_{\rm s}(\rho,\theta)$ are
respectively the real and imaginary parts of $w(z)$. Since the
coefficients $a_{k}$ are in fact independent of $\rho$ in this last set of
expressions, for $0<\rho<1$, and by construction have the same values as
those given by the previous set of expressions, in
Equation~(\ref{EQfrcoefs}), which define them as Fourier coefficients of
$f(\theta)$, we see that they are continuous from within as functions of
$\rho$ in the limit $\rho\to 1$, which we may then take.

According to the results of the previous section, we conclude therefore
that the two FC Fourier series we started with in this section are in fact
the DP Fourier series of the functions $f_{\rm c}(\rho,\theta)$ and
$f_{\rm s}(\rho,\theta)$ in the $\rho\to 1$ limit. Since the coefficients
uniquely identify the function almost everywhere, as discussed in the
introduction, we may now identify these two functions in the $\rho\to 1$
limit with the two functions $f_{\rm c}(\theta)$ and $f_{\rm s}(\theta)$
from which the coefficients $a_{k}$ were obtained in the first place, on
the unit circle. We conclude that the limits

\noindent
\begin{eqnarray*}
  f_{\rm c}(\theta)
  & = &
  \lim_{\rho\to 1}
  f_{\rm c}(\rho,\theta),
  \\
  f_{\rm s}(\theta)
  & = &
  \lim_{\rho\to 1}
  f_{\rm s}(\rho,\theta),
\end{eqnarray*}

\noindent
exist and hold almost everywhere over the unit circle. We may consider the
$\rho\to 1$ limits of the functions $f_{\rm c}(\rho,\theta)$ and $f_{\rm
  s}(\rho,\theta)$ to be the maximally smooth members of sets of functions
that are zero-measure equivalent and lead to the same set of Fourier
coefficients. The functions $f_{\rm c}(\rho,\theta)$ and $f_{\rm
  s}(\rho,\theta)$, which according to our definitions are the Fourier
Conjugate functions of each other, are also known as Harmonic Conjugate
functions~\cite{HarmConj}, since they are the real and imaginary parts of
an analytic function, and hence are both harmonic functions on the real
plane. They can be obtained from one another by the Hilbert
transform~\cite{HilbTran}. In the $\rho\to 1$ limit we may also get
functions $f_{\rm c}(\theta)$ and $f_{\rm s}(\theta)$ which are
restriction to the unit circle of Harmonic Conjugate functions, if the
function $w(z)$ is analytic at the limiting point. But in any case they
are Fourier Conjugate to each other.

Up to this point we have established that every given pair of FC Fourier
series, obtained from a given pair of FC real functions, such that both
converge together on at least one point on the unit circle, corresponds to
a specific inner analytic function, whose Taylor series converges on at
least that point in the unit circle, and which reproduces the original
pair of FC real function almost everywhere when one takes the $\rho\to 1$
limit from within the open unit disk to the unit circle.

Furthermore, we also see that we may work this argument in reverse. In
other words, given an arbitrary inner analytic function $w(z)$ whose
Taylor series converges on at least one point in the unit circle and which
is well-defined almost everywhere over that circle, we may construct from
it the DP Fourier series of two related functions. These two DP Fourier
series are just the real and imaginary parts of the Taylor series of the
analytic function $w(z)$, in the $\rho\to 1$ limit. Finally, since the
coefficients $a_{k}$ are continuous functions of $\rho$ from within at the
unit disk, and since the $S_{z}$ series converges on at least one point on
the unit circle, we may also conclude that the two corresponding FC
Fourier series converge together on that point of the unit circle.

This completes the establishment of a one-to-one correspondence between,
on the one hand, pairs of real FC Fourier series that converge together on
at least a single point of the periodic interval and, on the other hand,
inner analytic functions that converge on at least one point of the unit
circle and are well-defined almost everywhere over that circle. Note,
however, that the analytic side of this correspondence is the more
powerful one, because given the Fourier coefficients $a_{k}$ we may be
able to define a convergent power series $S_{z}$ and thus an inner
analytic function $w(z)$ even if the corresponding Fourier series diverge
everywhere on the periodic interval.

This correspondence can be a useful tool, as it may make it easier to
determine the convergence or lack thereof of given Fourier series. It may
also be used to recover from the coefficients $a_{k}$ the function which
generates a given Fourier series, even if that series is divergent. The
way to do this is simply to determine the corresponding inner analytic
function and then calculate the $\rho\to 1$ limit of the real and
imaginary parts of that function. In Appendix~\ref{APPexamplim} we will
give a few simple examples of this type of procedure.

\section{The Orthogonality Relations}

The two most central elements of the structure of Fourier theory are the
set of orthogonality relations and the completeness relation. Let us then
show that these also follow from the structure of complex analysis. We
start with the orthogonality relations, which we already gave in
Equation~(\ref{EQorthrels}) of the introduction. Of course the integrals
involved are simple ones, and can be calculated by elementary means. Our
objective here, however, is not to just calculate them but to show that
they are a consequence of the analytic structure of the complex plane. We
can do this by simply considering the Cauchy integral formulas for the
coefficients $a_{k'}$ of the Taylor series of a simple power $z^{k}$, with
$k\geq 0$,

\noindent
\begin{eqnarray*}
  a_{k'}
  & = &
  \frac{1}{k!}\,
  \frac{d^{k'}z^{k}}{dz^{k'}}(0)
  \\
  & = &
  \frac{1}{2\pi\ii}
  \oint_{C}dz\,
  \frac{z^{k}}{z^{k'+1}},
\end{eqnarray*}

\noindent
where $C$ is a circle or radius $\rho$ centered at the origin, with
$0<\rho\leq 1$. On the one hand, if $k'>k$ we get $a_{k'}=0$ due to the
multiple differentiation of the power function, which is differentiated
more times than the power itself. On the other hand, if $k'<k$ we get
$a_{k'}=0$ when we calculate the derivatives and apply the result at zero,
since in this case there is always at least one factor of $z$ left, or
alternatively due to the Cauchy-Goursat theorem, because in this case the
integrand is analytic and thus the integral is zero. If $k'=k$, however,
we get $a_{k'}=1$, which we can get either directly from the result of the
differentiation, or from the fact that in this case the integral is given
by

\begin{displaymath}
  \oint_{C}dz\,
  \frac{1}{z}
  =
  2\pi\ii,
\end{displaymath}

\noindent
as one can easily verify, either directly or by the residues theorem. In
any case we get the result

\begin{displaymath}
  \frac{1}{2\pi\ii}
  \oint_{C}dz\,
  \frac{z^{k}}{z^{k'+1}}
  =
  \delta_{kk'}.
\end{displaymath}

\noindent
We now write the integral explicitly on the circle of radius $\rho$, using
$\theta$ as the integration variable, with $z=\rho\exp(\ii\theta)$ and
thus with $dz=\ii zd\theta$,

\noindent
\begin{eqnarray*}
  \delta_{kk'}
  & = &
  \frac{1}{2\pi\ii}
  \int_{-\pi}^{\pi}d\theta\,
  \ii\,
  \rho^{k}
  \e{\ii k\theta}
  \rho^{-k'}
  \e{-\ii k'\theta}
  \\
  & = &
  \frac{\rho^{k-k'}}{2\pi}
  \int_{-\pi}^{\pi}d\theta\,
  \left[
    \cos(k\theta)
    +
    \ii
    \sin(k\theta)
  \right]
  \left[
    \cos(k'\theta)
    -
    \ii
    \sin(k'\theta)
  \right]
  \\
  & = &
  \frac{\rho^{k-k'}}{2\pi}
  \left\{
    \left[
      \int_{-\pi}^{\pi}d\theta\,
      \cos(k\theta)
      \cos(k'\theta)
      +
      \int_{-\pi}^{\pi}d\theta\,
      \sin(k\theta)
      \sin(k'\theta)
    \right]
  \right.
  +
  \\
  &   &
  \hspace{3.5em}
  +
  \left.
    \ii
    \left[
      \int_{-\pi}^{\pi}d\theta\,
      \sin(k\theta)
      \cos(k'\theta)
      -
      \int_{-\pi}^{\pi}d\theta\,
      \cos(k\theta)
      \sin(k'\theta)
    \right]
  \right\}.
\end{eqnarray*}

\noindent
One can see that the integrals in the imaginary part are zero due to
parity arguments. In fact, these constitute some of the orthogonality
relations, those including sines and cosines. We are left with

\begin{equation}\label{EQorthog}
  \delta_{kk'}
  =
  \frac{\rho^{k-k'}}{2\pi}
  \left[
    \int_{-\pi}^{\pi}d\theta\,
    \cos(k\theta)
    \cos(k'\theta)
    +
    \int_{-\pi}^{\pi}d\theta\,
    \sin(k\theta)
    \sin(k'\theta)
  \right].
\end{equation}

\noindent
For $k=0=k'$ the second term vanishes, and the equation becomes a simple
identity, which is in fact one of the other orthogonality relations. If,
on the other hand, we have $k+k'>0$, we now consider the integral, on the
same circuit,

\begin{displaymath}
  \oint_{C}dz\,
  z^{k}z^{k'-1}
  =
  0,
\end{displaymath}

\noindent
which is zero due to the Cauchy-Goursat theorem, since the integrand is
analytic within the circle $C$ for $k+k'\geq 1$. Writing the integral
explicitly on the circle we get

\noindent
\begin{eqnarray*}
  0
  & = &
  \int_{-\pi}^{\pi}d\theta\,
  \ii\,
  \rho^{k}
  \e{\ii k\theta}
  \rho^{k'}
  \e{\ii k'\theta}
  \\
  & = &
  \ii\rho^{k+k'}
  \int_{-\pi}^{\pi}d\theta\,
  \left[
    \cos(k\theta)
    +
    \ii
    \sin(k\theta)
  \right]
  \left[
    \cos(k'\theta)
    -
    \ii
    \sin(k'\theta)
  \right]
  \\
  & = &
  \ii\rho^{k+k'}
  \left\{
    \left[
      \int_{-\pi}^{\pi}d\theta\,
      \cos(k\theta)
      \cos(k'\theta)
      -
      \int_{-\pi}^{\pi}d\theta\,
      \sin(k\theta)
      \sin(k'\theta)
    \right]
  \right.
  +
  \\
  &   &
  \hspace{3.5em}
  +
  \left.
    \ii
    \left[
      \int_{-\pi}^{\pi}d\theta\,
      \sin(k\theta)
      \cos(k'\theta)
      +
      \int_{-\pi}^{\pi}d\theta\,
      \cos(k\theta)
      \sin(k'\theta)
    \right]
  \right\}.
\end{eqnarray*}

\noindent
Once again the integrals in the imaginary part are zero by parity
arguments, and thus we are left with

\begin{displaymath}
  \int_{-\pi}^{\pi}d\theta\,
  \cos(k\theta)
  \cos(k'\theta)
  =
  \int_{-\pi}^{\pi}d\theta\,
  \sin(k\theta)
  \sin(k'\theta),
\end{displaymath}

\noindent
since $\rho\neq 0$. If we now go back to our previous expression in
Equation~(\ref{EQorthog}) we see that the two integrals that appear there
are equal to each other, so we may write that

\noindent
\begin{eqnarray*}
  \frac{1}{\pi}
  \int_{-\pi}^{\pi}d\theta\,
  \cos(k\theta)
  \cos(k'\theta)
  & = &
  \rho^{k'-k}
  \delta_{kk'}
  \\
  & = &
  \delta_{kk'},
  \\
  \frac{1}{\pi}
  \int_{-\pi}^{\pi}d\theta\,
  \sin(k\theta)
  \sin(k'\theta)
  & = &
  \rho^{k'-k}
  \delta_{kk'}
  \\
  & = &
  \delta_{kk'},
\end{eqnarray*}

\noindent
where the factor involving $\rho$ is irrelevant since the right-hand sides
are only different from zero if $k=k'$. We get therefore the complete set
of orthogonality relations

\noindent
\begin{eqnarray*}
  \frac{1}{\pi}
  \int_{-\pi}^{\pi}d\theta\,
  \cos(k\theta)
  \cos(k'\theta)
  & = &
  \delta_{kk'},
  \\
  \frac{1}{\pi}
  \int_{-\pi}^{\pi}d\theta\,
  \sin(k\theta)
  \sin(k'\theta)
  & = &
  \delta_{kk'},
  \\
  \frac{1}{\pi}
  \int_{-\pi}^{\pi}d\theta\,
  \cos(k\theta)
  \sin(k'\theta)
  & = &
  0,
\end{eqnarray*}

\noindent
for $k\geq 1$ and $k'\geq 1$, which are the relevant values for DP Fourier
series. It is interesting to note that the orthogonality relations are
valid on the circle of radius $\rho$ with $0<\rho\leq 1$, without the need
to actually take the $\rho\to 1$ limit. We thus get a bit more than we
bargained for in this case, for it would have been sufficient to establish
these relation only on the unit circle. Note that this is different from
what happened during the calculation of the coefficients of the Fourier
series. However, in both cases the results come from the Cauchy integral
formulas and the Cauchy-Goursat theorem, and in either case the same real
integral appears, defining the usual scalar product in the space of real
functions on the periodic interval.

\section{The Completeness Relation}\label{SECcomprel}

Let us now show that the completeness relation of the Fourier basis also
follows from the structure of complex analysis. In order to do this, we
must first show that the Dirac delta ``function'' can be represented in
terms of the analytic structure within the open unit disk. We denote the
Dirac delta ``function'' centered at $\theta=\theta_{1}$ on the unit
circle by $\delta(\theta-\theta_{1})$. The definition of this mathematical
object is that it is a symbolic representation of a limiting process which
has the following four properties:

\begin{enumerate}

\item $\delta(\theta-\theta_{1})$ tends to zero when one takes the
  defining limit with $\theta\neq\theta_{1}$;

\item $\delta(\theta-\theta_{1})$ diverges to positive infinity when one
  takes the defining limit with $\theta=\theta_{1}$;

\item in the defining limit the integral

  \begin{displaymath}
    \int_{a}^{b}d\theta\,
    \delta(\theta-\theta_{1})
    =
    1,
  \end{displaymath}

  \noindent
  has the value shown, for any interval $(a,b)$ which contains the point
  $\theta_{1}$;

\item given any continuous function $g(\theta)$, in the defining limit the
  integral

  \begin{displaymath}
    \int_{a}^{b}d\theta\,
    g(\theta)
    \delta(\theta-\theta_{1})
    =
    g(\theta_{1}),
  \end{displaymath}

  \noindent
  has the value shown, for any interval $(a,b)$ which contains the point
  $\theta_{1}$.

\end{enumerate}

\noindent
Of course no real function exists that can have all these properties,
which justifies the quotes in which we wrap the word ``function'' when
referring to it. In order to construct the Dirac delta ``function'' we
must first give an object or set of objects over which the limiting
process can be defined, and then define that limiting process. In order to
fulfill this program, we consider the complex function given within the
open unit disk by

\begin{displaymath}
  w_{\delta}(z)
  =
  \frac{1}{2\pi}
  -
  \frac{1}{\pi}\,
  \frac{z}{z-z_{1}},
\end{displaymath}

\noindent
as well as its restrictions to circles of radius $\rho$ centered at the
origin, with $0<\rho<1$, where $z=\rho\exp(\ii\theta)$ and where
$z_{1}=\exp(\ii\theta_{1})$ is a point on the unit circle. This function
is analytic within the open unit disk, but it is not an inner analytic
function, because $w_{\delta}(0)$ is not zero. However, we may write it in
terms of another function $w(z)$ as

\noindent
\begin{eqnarray*}
  w_{\delta}(z)
  & = &
  \frac{1}{2\pi}
  +
  w(z),
  \\
  w(z)
  & = &
  -\,
  \frac{1}{\pi}\,
  \frac{z}{z-z_{1}}.
\end{eqnarray*}

\noindent
Strictly speaking, $w(z)$ is not an inner analytic function either,
because it does not reduce to a real function over the real axis. However,
it does reduce to a real function over the straight line $z=\chi z_{1}$,
with real $\chi$, since in this case we have

\begin{displaymath}
  w(z)
  =
  -\,
  \frac{1}{\pi}\,
  \frac{\chi}{\chi-1}.
\end{displaymath}

\noindent
We see therefore that $w(z)$ is an inner analytic function rotated around
the origin by the angle $\theta_{1}$ associated to $z_{1}$. Therefore this
is just a simple extension of the structure we defined here. The limiting
process to be used for the definition is just the $\rho\to 1$ limit to the
unit circle. We will now show that the real part of $w_{\delta}(z)$, taken
on the $\rho\to 1$ limit, satisfies all the required properties defining
the Dirac delta ``function''. In order to recover the real and imaginary
parts of this complex function, we must now rationalize it,

\noindent
\begin{eqnarray*}
  w_{\delta}(z)
  & = &
  \frac{1}{2\pi}
  -
  \frac{1}{\pi}\,
  \frac{z(z^{*}-z_{1}^{*})}{(z-z_{1})(z^{*}-z_{1}^{*})}
  \\
  & = &
  \frac{1}{2\pi}
  -
  \frac{1}{\pi}\,
  \frac
  {\left(\rho^{2}-zz_{1}^{*}\right)}
  {\rho^{2}-(zz_{1}^{*}+z^{*}z_{1})+1}
  \\
  & = &
  \frac{1}{2\pi}
  -
  \frac{1}{\pi}\,
  \frac
  {\rho^{2}-\rho\cos(\Delta\theta)-\ii\rho\sin(\Delta\theta)}
  {\rho^{2}-2\rho\cos(\Delta\theta)+1},
\end{eqnarray*}

\noindent
where $\Delta\theta=\theta-\theta_{1}$. We now examine the real part of
this function,

\begin{displaymath}
  \Re[w_{\delta}(z)]
  =
  \frac{1}{2\pi}
  -
  \frac{1}{\pi}\,
  \frac
  {\rho\left[\rho-\cos(\Delta\theta)\right]}
  {\left(\rho^{2}+1\right)-2\rho\cos(\Delta\theta)}.
\end{displaymath}

\noindent
If we now take the limit $\rho\to 1$, under the assumption that
$\Delta\theta\neq 0$, we get

\noindent
\begin{eqnarray*}
  \lim_{\rho\to 1}\Re[w_{\delta}(z)]
  & = &
  \frac{1}{2\pi}
  -
  \frac{1}{\pi}\,
  \frac
  {1-\cos(\Delta\theta)}
  {2-2\cos(\Delta\theta)}
  \\
  & = &
  \frac{1}{2\pi}
  -
  \frac{1}{2\pi}
  \\
  & = &
  0,
\end{eqnarray*}

\noindent
which is the correct value for the case of the Dirac delta ``function''.
Thus we see that the first property holds.

If, on the other hand, we calculate $\Re[w_{\delta}(z)]$ for
$\Delta\theta=0$ and $\rho<1$ we obtain

\noindent
\begin{eqnarray*}
  \Re[w_{\delta}(z)]
  & = &
  \frac{1}{2\pi}
  -
  \frac{1}{\pi}\,
  \frac
  {\rho(\rho-1)}
  {(\rho-1)^{2}}
  \\
  & = &
  \frac{1}{2\pi}
  -
  \frac{1}{\pi}\,
  \frac{\rho}{\rho-1},
\end{eqnarray*}

\noindent
which diverges to positive infinity as $\rho\to 1$ from below, as it
should in order to represent the singular Dirac delta ``function''. This
establishes that the second property holds.

We then calculate the integral of $\Re[w_{\delta}(z)]$ over the circle of
radius $\rho<1$, which is given by

\noindent
\begin{eqnarray*}
  I_{1}
  & = &
  \frac{1}{2\pi}
  \int_{-\pi}^{\pi}d\theta\,
  \rho
  \left\{
    1
    -
    \frac
    {2\rho\left[\rho-\cos(\Delta\theta)\right]}
    {\left(\rho^{2}+1\right)-2\rho\cos(\Delta\theta)}
  \right\}
  \\
  & = &
  \frac{\rho}{2\pi}
  \int_{-\pi}^{\pi}d\theta\,
  \frac
  {\left(1-\rho^{2}\right)}
  {\left(\rho^{2}+1\right)-2\rho\cos(\Delta\theta)}
  \\
  & = &
  \frac{\left(1-\rho^{2}\right)}{4\pi}
  \int_{-\pi}^{\pi}d\theta\,
  \frac{1}{\left[\left(\rho^{2}+1\right)/(2\rho)\right]-\cos(\Delta\theta)},
\end{eqnarray*}

\noindent
since $d(\Delta\theta)=d\theta$. Note that this is {\em not} the integral
of an analytic function over a closed contour, but the integral of a real
function over the circle of radius $\rho$. This real integral over
$\theta$ can be calculated by residues. We introduce an auxiliary complex
variable $\xi=\lambda\exp(\ii\Delta\theta)$, which becomes simply
$\exp(\ii\Delta\theta)$ on the unit circle $\lambda=1$. We have
$d\xi=\ii\xi d\theta$, and so we may write the integral as

\noindent
\begin{eqnarray*}
  \int_{-\pi}^{\pi}d\theta\,
  \frac{1}{\left[\left(1+\rho^{2}\right)/(2\rho)\right]-\cos(\Delta\theta)}
  & = &
  \oint_{C}d\xi\,
  \frac{1}{\ii\xi}\,
  \frac{2}{\left[\left(1+\rho^{2}\right)/\rho\right]-\xi-1/\xi}
  \\
  & = &
  2\ii
  \oint_{C}d\xi\,
  \frac{1}{1-\left[\left(1+\rho^{2}\right)/\rho\right]\xi+\xi^{2}},
\end{eqnarray*}

\noindent
where the integral is now over the unit circle $C$ in the complex $\xi$
plane. The two roots of the quadratic polynomial on $\xi$ in the
denominator are given by

\noindent
\begin{eqnarray*}
  \xi_{+}
  & = &
  1/\rho,
  \\
  \xi_{-}
  & = &
  \rho.
\end{eqnarray*}

\noindent
Since $\rho<1$, only the pole corresponding to $\xi_{-}$ lies inside the
integration contour, so we get for the integral

\noindent
\begin{eqnarray*}
  \int_{-\pi}^{\pi}d\theta\,
  \frac{1}{\left[\left(1+\rho^{2}\right)/(2\rho)\right]-\cos(\Delta\theta)}
  & = &
  2\ii(2\pi\ii)
  \lim_{\xi\to\rho}
  \frac{1}{\xi-1/\rho}
  \\
  & = &
  4\pi\,
  \frac{\rho}{\left(1-\rho^{2}\right)}.
\end{eqnarray*}

\noindent
It follows that we have for the integral $I_{1}$

\noindent
\begin{eqnarray*}
  I_{1}
  & = &
  \frac{\left(1-\rho^{2}\right)}{4\pi\rho}\,
  4\pi\,
  \frac{\rho}{\left(1-\rho^{2}\right)}
  \\
  & = &
  1,
\end{eqnarray*}

\noindent
independently of $\rho$, including therefore the $\rho\to 1$ limit. Once
we have this result, and since the integrand goes to zero everywhere on
the unit circle except at $\Delta\theta=0$, which means that
$\theta=\theta_{1}$, the integral can be changed to one over any open
interval on the unit circle containing the point $\theta_{1}$, without any
change in its limiting value. This establishes that the third property
holds.

In order to establish the fourth and last property, we take an essentially
arbitrary inner analytic function $\gamma(z)$, with the single additional
restriction that it be well-defined at the point $z_{1}$, in the sense
that its $\rho\to 1$ limit exists at $z_{1}$. This inner analytic function
corresponds to a pair of FC real functions on the unit circle, both of
which are well-defined at $z_{1}$. We now consider the following integral
over the circle of radius $\rho<1$,

\noindent
\begin{eqnarray*}
  I_{2}
  & = &
  \frac{1}{2\pi}
  \int_{-\pi}^{\pi}d\theta\,
  \rho\,
  \gamma(z)
  \left\{
    1
    -
    \frac
    {2\rho\left[\rho-\cos(\Delta\theta)\right]}
    {\left(\rho^{2}+1\right)-2\rho\cos(\Delta\theta)}
  \right\}
  \\
  & = &
  \frac{\rho}{2\pi}
  \int_{-\pi}^{\pi}d\theta\,
  \gamma(z)\,
  \frac
  {\left(1-\rho^{2}\right)}
  {\left(\rho^{2}+1\right)-2\rho\cos(\Delta\theta)}
  \\
  & = &
  \frac{\left(1-\rho^{2}\right)}{4\pi}
  \int_{-\pi}^{\pi}d\theta\,
  \frac
  {\gamma(z)}
  {\left[\left(\rho^{2}+1\right)/(2\rho)\right]-\cos(\Delta\theta)},
\end{eqnarray*}

\noindent
since $d(\Delta\theta)=d\theta$. Note once more that this is {\em not} the
integral of an analytic function over a closed contour, but two integrals
of real functions, given by the real and imaginary parts of $\gamma(z)$,
over the circle of radius $\rho$. These real integrals over $\theta$ can
be calculated by residues, exactly like the one which appeared before in
the case of $I_{1}$. The calculation is exactly the same except for the
extra factor of $\gamma(z)$ to be taken into consideration when
calculating the residue, so that we may write directly that

\noindent
\begin{eqnarray*}
  \int_{-\pi}^{\pi}d\theta\,
  \frac
  {\gamma(z)}
  {\left[\left(\rho^{2}+1\right)/(2\rho)\right]-\cos(\Delta\theta)}
  & = &
  2\ii(2\pi\ii)
  \lim_{\xi\to\rho}
  \frac{\gamma(z)}{\xi-1/\rho}
  \\
  & = &
  4\pi\,
  \frac{\rho}{\left(1-\rho^{2}\right)}
  \lim_{\xi\to\rho}\gamma(z).
\end{eqnarray*}

\noindent
Note now that since $\xi=\lambda\exp(\ii\Delta\theta)$ and we must take
the limit $\xi\to\rho$, we in fact have that in that limit

\begin{displaymath}
  \lambda\e{\ii\Delta\theta}
  =
  \rho,
\end{displaymath}

\noindent
which implies that $\lambda=\rho$ and that $\Delta\theta=0$. We must
therefore write $\gamma(z)$ at the point given by $\rho$ and
$\Delta\theta=0$, that is, at the point given by $\rho$ and $\theta_{1}$,

\begin{displaymath}
  \int_{-\pi}^{\pi}d\theta\,
  \frac
  {\gamma(z)}
  {\left[\left(\rho^{2}+1\right)/(2\rho)\right]-\cos(\Delta\theta)}
  =
  4\pi\,
  \frac{\rho}{\left(1-\rho^{2}\right)}\,
  \gamma(\rho,\theta_{1}).
\end{displaymath}

\noindent
It follows that we have for the integral $I_{2}$

\noindent
\begin{eqnarray*}
  I_{2}
  & = &
  \frac{\left(1-\rho^{2}\right)}{4\pi\rho}\,
  4\pi\,
  \frac{\rho}{\left(1-\rho^{2}\right)}\,
  \gamma(\rho,\theta_{1})
  \\
  & = &
  \gamma(\rho,\theta_{1}).
\end{eqnarray*}

\noindent
Finally, we may now take the $\rho\to 1$ limit, since
$\gamma(\rho,\theta)$ is well-defined in that limit, and thus obtain

\begin{displaymath}
  \lim_{\rho\to 1}
  I_{2}
  =
  \gamma(1,\theta_{1}).
\end{displaymath}

\noindent
Once we have this result, and since the integrand goes to zero everywhere
on the unit circle except at $\Delta\theta=0$, which means that
$\theta=\theta_{1}$, the integral can be changed to one over any open
interval on the unit circle containing the point $\theta_{1}$, without any
change in its value. This establishes that the fourth and last property
holds. We may then write symbolically that

\begin{displaymath}
  \delta(\theta-\theta_{1})
  =
  \lim_{\rho\to 1}
  \Re[w_{\delta}(z)].
\end{displaymath}

\noindent
Note that in order to obtain this result it was not necessary to assume
that $\gamma(z)$ is continuous at $z_{1}$ in the direction of $\theta$
along the unit circle. It was necessary to assume only that $\gamma(z)$ is
continuous as a function of $\rho$, in the direction perpendicular to the
unit circle. One can see therefore that, once more, we get a bit more than
we bargained for, because we were able to establish the result with
slightly weaker hypotheses than at first expected.

We are now in a position to establish the completeness relation using this
representation of the Dirac delta ``function''. If we use once again the
Cauchy integral formulas for $w(z)$ we get for the coefficients of the
Taylor expansion of $w(z)$

\noindent
\begin{eqnarray*}
  a_{k}
  & = &
  \frac{1}{2\pi\ii}
  \oint_{C}dz\,
  \frac{w(z)}{z^{k+1}}
  \\
  & = &
  \frac{1}{2\pi^{2}\ii}
  \oint_{C}dz\,
  \frac{1}{z^{k}}\,
  \frac{(-1)}{z-z_{1}},
\end{eqnarray*}

\noindent
for $k\geq 1$, since $w(z)$ is a rotated inner analytic function. We now
observe that the second ratio in the integrand can be understood as the
sum of a geometric series, which is convergent so long as $\rho<1$,

\begin{displaymath}
  \frac{(-1)}{z-z_{1}}
  =
  \frac{1}{z_{1}}
  \sum_{n=0}^{\infty}
  \left(\frac{z}{z_{1}}\right)^{n},
\end{displaymath}

\noindent
so that we may now write

\noindent
\begin{eqnarray*}
  a_{k}
  & = &
  \frac{1}{2\pi^{2}\ii}
  \oint_{C}dz\,
  \frac{1}{z^{k}}\,
  \frac{1}{z_{1}}
  \sum_{n=0}^{\infty}
  \left(\frac{z}{z_{1}}\right)^{n}
  \\
  & = &
  \frac{1}{2\pi^{2}\ii}
  \sum_{n=0}^{\infty}
  \frac{1}{z_{1}^{n+1}}
  \oint_{C}dz\,
  z^{n-k},
\end{eqnarray*}

\noindent
since a convergent power series can always be integrated term-by-term. As
we have already discussed before, in the previous section, the remaining
integral is zero except if $n=k-1$, in which case it has the value
$2\pi\ii$. Note that this condition relating $n$ and $k$ can always be
satisfied since $k\geq 1$. We therefore get for the coefficients

\noindent
\begin{eqnarray*}
  a_{k}
  & = &
  \frac{1}{2\pi^{2}\ii}
  \sum_{n=0}^{\infty}
  \frac{1}{z_{1}^{n+1}}\,
  2\pi\ii\,\delta_{n,k-1}
  \\
  & = &
  \frac{1}{\pi z_{1}^{k}}.
\end{eqnarray*}

\noindent
As a result, we get for the Taylor expansion of $w(z)$

\begin{displaymath}
  w(z)
  =
  \frac{1}{\pi}
  \sum_{k=1}^{\infty}
  \left(\frac{z}{z_{1}}\right)^{k}.
\end{displaymath}

\noindent
We now write both $z$ and $z_{1}$ in polar coordinates, to obtain

\noindent
\begin{eqnarray*}
  w(z)
  & = &
  \frac{1}{\pi}
  \sum_{k=1}^{\infty}
  \rho^{k}
  \e{\ii k\theta}
  \e{-\ii k\theta_{1}}
  \\
  & = &
  \frac{1}{\pi}
  \sum_{k=1}^{\infty}
  \rho^{k}
  \left[
    \cos(k\theta)
    +
    \ii
    \sin(k\theta)
  \right]
  \left[
    \cos(k\theta_{1})
    -
    \ii
    \sin(k\theta_{1})
  \right]
  \\
  & = &
  \frac{1}{\pi}
  \sum_{k=1}^{\infty}
  \rho^{k}
  \left\{
    \rule{0em}{3ex}
    \left[
      \cos(k\theta)
      \cos(k\theta_{1})
      +
      \sin(k\theta)
      \sin(k\theta_{1})
    \right]
  \right.
  +
  \\
  &   &
  \hspace{4.5em}
  +
  \ii
  \left.
    \left[
      \sin(k\theta)
      \cos(k\theta_{1})
      -
      \cos(k\theta)
      \sin(k\theta_{1})
    \right]
    \rule{0em}{3ex}
  \right\}.
\end{eqnarray*}

\noindent
If we now write the real part of $w_{\delta}(z)$ we get

\begin{displaymath}
  \Re[w_{\delta}(z)]
  =
  \frac{1}{2\pi}
  +
  \frac{1}{\pi}
  \sum_{k=1}^{\infty}
  \rho^{k}
  \left[
    \cos(k\theta)
    \cos(k\theta_{1})
    +
    \sin(k\theta)
    \sin(k\theta_{1})
  \right],
\end{displaymath}

\noindent
and, if we then take the $\rho\to 1$ limit, we get the expression

\begin{displaymath}
  \delta(\theta-\theta_{1})
  =
  \frac{1}{2\pi}
  +
  \frac{1}{\pi}
  \sum_{k=1}^{\infty}
  \left[
    \cos(k\theta)
    \cos(k\theta_{1})
    +
    \sin(k\theta)
    \sin(k\theta_{1})
  \right],
\end{displaymath}

\noindent
which is the completeness relation in its usual form, a bilinear form on
the Fourier basis functions, at two separate points along the unit circle.
Note that the constant function, which is an element of the complete
Fourier basis, is included in the first term. Note also that this time it
was necessary to take the $\rho\to 1$ limit, and that the completeness of
the Fourier basis is valid only on the unit circle. This is to be
expected, of course, since the unit circle is where the corresponding
space of real functions, which is generated by the basis, is defined.

\section{Limits from Within}

We have established that every DP Fourier series that converges on one
point together with its FC series corresponds to an inner analytic
function. Whenever it is possible to take the limit to the unit circle of
restrictions of this analytic function to circles of smaller radii,
centered at the origin, it gives us back the real function that
corresponds to the coefficients of the series, even if the Fourier series
itself is not convergent. We may also interpret this type of limit as a
collection of point-by-point limits taken in the radial direction to each
point of the unit circle. Let us discuss now under what conditions we may
take such limits and what we can learn from them.

Our ability to take the limits to the unit circle depends on whether or
not the inner analytic function $w(z)$ is well-defined on the unit circle,
and therefore on whether or not it has singularities on that circle, as
well as on the nature of these singularities. If there are no
singularities at all on the unit circle, then $w(z)$ is analytic over the
whole unit circle and therefore continuous there. In this case it is
always possible to take the limit, to all points of the circle, and they
will always result in a pair of $C^{\infty}$ real functions on the unit
circle. Also, in this case it is true that the Taylor series of $w(z)$ and
the corresponding pair of FC series converge everywhere on the unit
circle. In fact, as we will see shortly, they all converge absolutely and
uniformly, and this is the situation which we characterize and refer to as
that of strong convergence. The $C^{\infty}$ functions thus obtained on
the unit circle are those that give the coefficients $a_{k}$.

Let us suppose now that there is a finite number of singularities on the
unit circle, which are therefore all isolated singularities. On the open
subsets of the unit circle between two adjacent singularities $w(z)$ is
still analytic, and hence continuous. Therefore, within these open subsets
the limits to the unit circle may still be taken, resulting in segments of
$C^{\infty}$ real functions, and reproducing almost everywhere over the
circle the real functions that originated the $a_{k}$ coefficients.
Therefore we learn that any Fourier series that converges to a sectionally
continuous and differentiable function, and is such that the corresponding
inner analytic function has a finite number of singularities on the unit
circle, in fact converges to a sectionally $C^{\infty}$ function, possibly
with an increased number of sections.

At the points of singularity, if the Fourier series converge, then by
Abel's theorem~\cite{AbelTheo} they converge to the limits of $w(z)$ at
these points, taken from within the open unit disk, so that in this case
it is also possible to take the limits. Since the real functions that
generate the coefficients by means of the real integrals are determined
uniquely only almost everywhere, if the series diverge at the singular
points we may still adopt the limits of $w(z)$ from within as the values
of the corresponding functions at those points, so long as these limits
exist and are finite. This will be the case if the singularities do not
involve divergences to infinity, so that the inner analytic functions are
still well-defined on their location, although they are not analytic
there. We will call such singularities, at which we may still take the
limits to the unit circle, {\em soft singularities}.

Usually one tends to think of singularities in analytic functions in terms
only of {\em hard singularities} such as poles, which always involve
divergences to infinity when one approaches them. The paradigmatic hard
singularities are poles such as

\begin{displaymath}
  \frac{1}{(z-z_{1})^{n}}
  =
  \frac{\e{-\ii n\alpha}}{R^{n}},
\end{displaymath}

\noindent
with $n\geq 1$, where $z-z_{1}=R\exp(\ii\alpha)$ with real $R$ and
$\alpha$, and where $z_{1}$ is a point on the unit circle. We will think
of the integer $n$ as the {\em degree of hardness} of the singularity.
Note that such singularities are not integrable along arbitrary lines
passing through the singular point. A somewhat less hard singularity is
the logarithmic one given by

\begin{displaymath}
  \ln(z-z_{1})
  =
  \ln(R)+\ii\alpha,
\end{displaymath}

\noindent
which is also hard but much less so that the poles, going to infinity much
slower. We will call this a {\em borderline hard} singularity, or a hard
singularity of degree zero. Note that this singularity is integrable along
arbitrary lines passing through the singular point. Soft singularities are
those for which the limit to the singular point still exists, for example
those such as

\begin{displaymath}
  (z-z_{1})\ln(z-z_{1}),
\end{displaymath}

\noindent
which we will call a {\em borderline soft} singularity, or a soft
singularity of degree zero. Progressively softer singularities can be
obtained by increasing $n$ in the generic soft singularity

\begin{displaymath}
  (z-z_{1})^{n+1}\ln(z-z_{1}),
\end{displaymath}

\noindent
where we should have $n\geq 1$. In this case we will think of the integer
$n$ as the {\em degree of softness} of the singularity. To complete this
classification of singularities, any essential singularity will be
classified as an infinitely hard one. A more precise and complete
definition of this classification of singularities will be given in the
follow-up paper~\cite{FTotCPII}.

If one has a countable infinity of singularities on the unit circle, but
they have only a finite number of accumulation points, then the situation
still does not change too much. In the open subsets of the unit circle
between every pair of consecutive singular points the limits still exist
and give us segments of $C^{\infty}$ functions. If the singularities are
all soft, then the limits exist everywhere over the unit circle. In any
case the limits give us back the real functions which generated the
coefficients, at least almost everywhere over the whole circle.

Even if we have a countably infinite set of singularities which are
distributed densely on the unit circle, so long as the singularities are
all soft the limits can still be taken everywhere. Of course, in this case
it is not to be expected that the resulting real functions will be
$C^{\infty}$ anywhere. If a given inner analytic function has the open
unit disk as the maximum disk of convergence of its Taylor series, and has
a dense set of hard singularities over the unit circle, then the limits do
not exist anywhere, and therefore the real functions are in fact not
defined at all by the limits to that circle. In this case, due to Abel's
theorem~\cite{AbelTheo}, the Taylor series must diverge everywhere over
the unit circle.

Since a borderline hard singularity such as $\ln(z-z_{1})$ can be obtained
by the differentiation of a borderline soft singularity such as
$(z-z_{1})\ln(z-z_{1})$, it is reasonable to presume that the transition
by the differentiation of $w(z)$ from a densely distributed countable
infinity of borderline soft singularities to a densely distributed
countable infinity of borderline hard singularities corresponds to real
functions that exist everywhere over the unit circle but that are not
differentiable anywhere.

The questions related to the issue of convergence will be discussed in
more detail in the follow-up paper~\cite{FTotCPII}, in which we will show
that the degrees of hardness or softness of the singularities are related
to the analytic character of the real functions, all the way from simple
continuity through $n$-fold or $C^{n}$ differentiability to $C^{\infty}$
status. In particular, borderline soft singularities are related to
everywhere continuous functions, so that the transition described above
should be related to real functions that are everywhere continuous but
nowhere differentiable. One such proposed function will be discussed
briefly in Appendix~\ref{APPriemannf}.

\section{Some Basic Convergence Results}

Let us now consider the convergence of DP Fourier series and of the
corresponding complex power series. In many cases the convergence of the
complex power series, which is the Taylor series of the inner analytic
function, and which is usually more easily established, will give us
information about the convergence of the DP Fourier series. However,
sometimes the reverse situation may also be realized. For organizational
reasons this discussion must be separated into several parts. We will
tackle in this paper only the extreme cases which we will qualify as those
of {\em strong divergence} and {\em strong convergence}, for which results
can be established in full generality. Further discussions of the
convergence issues will be presented in the follow-up
paper~\cite{FTotCPII}.

\subsection{Strong Divergence}

In order to discuss this case we must go all the way back to the case of
plain DP trigonometric series, without any consideration of whether or not
they are Fourier series. Hence, let us suppose that we have an arbitrary
pair of FC trigonometric series,

\noindent
\begin{eqnarray*}
  S_{\rm c}
  & = &
  \sum_{k=1}^{\infty}
  a_{k}\cos(k\theta),
  \\
  S_{\rm s}
  & = &
  \sum_{k=1}^{\infty}
  a_{k}\sin(k\theta).
\end{eqnarray*}

\noindent
Without any further assumptions, we may then construct the complex series
$S_{v}$ and $S_{z}$,

\noindent
\begin{eqnarray*}
  S_{v}
  & = &
  \sum_{k=1}^{\infty}
  a_{k}\e{\ii k\theta},
  \\
  S_{z}
  & = &
  \sum_{k=1}^{\infty}
  a_{k}z^{k}.
\end{eqnarray*}

\noindent
Let us suppose that $S_{z}$ fails to converge at a single point $z_{1}$
located strictly {\em inside} the open unit disk. Then it follows from the
basic convergence theorem that it cannot converge at any point outside the
circle centered at zero with its boundary going through the point $z_{1}$.
This is so because, if the series did converge at some point outside this
disk, them by that theorem it would converge everywhere in a larger disk,
which contains the point of divergence $z_{1}$. This is absurd, and
therefore we conclude that the series must be divergent at all points
strictly outside the closed disk defined by the point of divergence
$z_{1}$.

It follows therefore that in this case the $S_{z}$ series is divergent
everywhere over the whole unit circle. The same is then true for $S_{v}$,
of course, and similar conclusions may be drawn for the associated DP
trigonometric series. Note that while the convergence of $S_{v}$ at a
given point implies the convergence of both $S_{\rm c}$ and $S_{\rm s}$ at
that point, the lack of convergence of $S_{v}$ implies only that at least
one of the two DP trigonometric series diverges. The other one might still
converge. We may immediately conclude that at least one of the two series
in the FC pair is divergent. But in fact, one can see that both must be
divergent almost everywhere, by the following argument.

If we consider the absolute value of the terms of $S_{z}$ at $z_{1}$, we
have $\left|a_{k}z_{1}^{k}\right|=|a_{k}|\rho_{1}^{k}$, where
$|z_{1}|=\rho_{1}$ and $\rho_{1}<1$. Since $S_{z}$ is divergent at
$z_{1}$, $|a_{k}|\rho_{1}^{k}$ cannot go to zero faster than $1/k$,
because a limit to zero as $1/k^{1+\varepsilon}$ with any strictly
positive $\varepsilon$ would be sufficient for the series of the absolute
values of the terms of $S_{z}$ to converge, since in this case it can be
bounded from above by a convergent asymptotic integral, as is shown in
Section~\ref{APPevalconv} of Appendix~\ref{APPproofs}. This would imply
that $S_{z}$ is absolutely convergent and hence convergent, which is
absurd. This means that we must have

\begin{displaymath}
  |a_{k}|
  >
  \frac{1}{k\rho_{1}^{k}},
\end{displaymath}

\noindent
for $k>k_{m}$ and some minimum value $k_{m}$ of $k$. In other words,
$|a_{k}|$ must in fact diverge to infinity exponentially with $k$, in fact
just about as fast as $1/\rho_{1}^{k}$, which goes to infinity as
$k\to\infty$ since $\rho_{1}<1$. Since these coefficients are those of
both $S_{\rm c}$ and $S_{\rm s}$, and a limit of the coefficients to zero
as $k\to\infty$ is a necessary condition for convergence, both these
series must diverge almost everywhere, that is, everywhere except possibly
for a few special points, such as the sine series at $\theta=0$.

Therefore, in this case there are no trigonometric series at all that
converge almost everywhere on the unit circle. We say that such series are
divergent almost everywhere. Besides, since any Fourier series is a
particular case of trigonometric series, in this case there are no
almost-everywhere convergent Fourier series as well. This is true whenever
the maximum disk of convergence of $S_{z}$ around $z=0$ is smaller than
the open unit disk, which also means that the function $w(z)$ to which
$S_{z}$ converges must have at least one singularity strictly within that
disk.

In conclusion, given any DP trigonometric series, be it a Fourier series
or not, if the power series built from its coefficients is divergent
anywhere strictly within the unit disk, then the trigonometric series is
divergent almost everywhere in the periodic interval. The same follows if
it can be verified that the function $w(z)$ has a singularity strictly
within the open unit disk. We refer to this situation as one of strong
divergence.

Note that in this case we cannot assert whether or not the trigonometric
series is a Fourier series. In this case it is not possible to take any
limit from the interior of the maximum disk of convergence to the unit
circle, which becomes disconnected from the analytic region of the $S_{z}$
series. In fact, it is an interesting question whether or not, in this
situation, any real functions can exist that give the coefficients of the
trigonometric series by means of the usual integrals. While a negative
answer to this question seems to be compelling, for now we must leave this
here as a simple conjecture.

\subsection{Strong Convergence}

The simplest case of convergence of the DP Fourier series is that in which
there is strong convergence to a $C^{\infty}$ function. Let us assume that
we have a pair of FC Fourier series and use their coefficients to
construct the corresponding $S_{v}$ and $S_{z}$ series. From the basic
convergence theorem it follows that, if $S_{z}$ converges at a single
point $z_{1}$ strictly {\em outside} the closed unit disk, then it is
strongly convergent over that whole disk. This means that $w(z)$ is
analytic everywhere on the unit circle, and therefore that it is
$C^{\infty}$ there. In fact, it means that the maximum disk of convergence
of the $S_{z}$ series is larger than the closed unit disk, and contains
it. In this case the inner analytic function $w(z)$ has no singularities
at all over the whole closed unit disk.

Since the series $S_{v}$ is a restriction of $S_{z}$ to the unit circle,
it then follows that $S_{v}$ and the corresponding DP Fourier series are
strongly convergent to $C^{\infty}$ functions over the whole unit circle,
and hence everywhere over the whole periodic interval. In fact, since this
is true for the power series $S_{z}$, in this case the DP Fourier series
can be differentiated term-by-term any number of times, and this operation
will always result in other equally convergent series. This is the
strongest form of convergence one can hope for. It is interesting to note,
in passing, that Fourier series displaying this type of convergence appear
in the solutions of essentially all boundary value problems involving the
diffusion equation in Cartesian coordinates.

Note that in this strong convergence case there is no need to make the
distinction between trigonometric series and Fourier series. This is so
because, if a trigonometric series converges almost everywhere on the unit
circle to a function $f(\theta)$, then the orthogonality of the Fourier
basis formed by the sets of cosine and sine functions suffices to show
that the coefficients of the series are given in terms of the usual
integrals involving $f(\theta)$. Therefore, under these circumstances
every trigonometric series is a Fourier series.

\subsection{A Strong-Convergence Test}

The situation in the strong convergence case suggests the following strong
convergence test for an arbitrary DP Fourier series. Given for example a
cosine series, and the corresponding FC series, one may write down the
pair of extended series

\noindent
\begin{eqnarray*}
  S'_{\rm c}
  & = &
  \sum_{k=1}^{\infty}
  a_{k}\rho^{k}\cos(k\theta),
  \\
  \bar{S}'_{\rm c}
  & = &
  \sum_{k=1}^{\infty}
  a_{k}\rho^{k}\sin(k\theta),
\end{eqnarray*}

\noindent
where $\rho\geq 0$ is real. If it can be established that there is a
single point, with $\rho>1$ and any value of $\theta$, such that these
series are both convergent at that point, then the original cosine series
converges strongly to a $C^{\infty}$ function over its whole domain, that
is, the whole periodic interval. The values $\theta=0$ and $\theta=\pm\pi$
are sensible choices for the test, because the sine series is always
convergent for these values, so that only the cosine series has to be
actually tested. The same test holds if we start with a sine series, of
course. Additionally, the corresponding FC series are also strongly
convergent everywhere on the periodic interval to the respective FC
functions.

This same test can also be used to determine whether the given series is
strongly divergent, if we revert the condition on $\rho$ and use the test
with $\rho<1$, looking now for a point of divergence. In this case it is
enough that any one of the two DP trigonometric series of the FC pair
diverge at a single point in the interior of the unit disk to ensure that
both diverge almost everywhere over the unit circle.

\subsection{A Modified Ratio Test}

Let us suppose that the $S_{z}$ series satisfies the ratio convergence
test at a point on the unit circle. That immediately implies, of course,
that the series converges strongly to an inner analytic function strictly
inside the unit disk. However, one can state more than just that, because
the test actually implies that the corresponding inner analytic function
is also analytic everywhere over the unit circle.

As demonstrated in Section~\ref{APPratitest} of Appendix~\ref{APPproofs},
if $S_{z}$ satisfies the ratio test at a point on the unit circle then,
due to the basic convergence theorem, the maximum disk of convergence of
the series $S_{z}$ extends beyond the unit circle, and contains it. It
follows therefore that the $S_{z}$ series converges strongly to an
analytic function on the whole closed unit disk. It now follows once again
that the $S_{v}$ series and the corresponding DP Fourier series $S_{\rm
  c}$ and $S_{\rm s}$ converge strongly to $C^{\infty}$ functions over the
whole unit circle, and hence everywhere over the whole periodic interval.

This constitutes in fact a modified ratio test that can be used for
arbitrary DP Fourier series. Given any cosine or sine series, if the
coefficients of the series satisfy the ratio test (and we mean here the
{\em coefficients} $a_{k}$, not the terms of the series), then the series
is strongly convergent to a $C^{\infty}$ function on the whole periodic
interval. The same is valid for its FC Fourier series, of course. Note
that, since the test is to be applied to the coefficients, in this case it
is not necessary to test the two FC series independently.

\section{Conclusions}

We have shown that there is a close and deep relationship between real
Fourier series and analytic functions in the unit disk centered at the
origin of the complex plane. In fact, there is a one-to-one correspondence
between pairs of FC Fourier series and the Taylor series of inner analytic
functions, so long as in both cases one has convergence at least on a
single point of the unit circle. This allows one to use the powerful and
extremely well-known machinery of complex analysis to deal with
trigonometric series in general, and Fourier series in particular.

One rather unexpected result of this interaction is the fact that the
well-known formulas for the Fourier coefficients are in fact a consequence
of the structure of the complex analysis of inner analytic functions.
Furthermore, so are all the other basic elements of Fourier theory, also
in a somewhat unexpected way. We may therefore conclude, from what we have
shown here, that the following rather remarkable statement is true:

\vspace{3ex}

\noindent
\parbox{\textwidth} {\bf The whole structure of the Fourier basis and of
  the associated theory is contained within the complex analysis of the
  set of inner analytic functions, including the orthogonality relations,
  the completeness relation, the form of the coefficients of the series
  and the form of the scalar product in the space of real functions
  generated by the basis.}

\vspace{3ex}

\noindent
One does not usually associate very weakly convergent Fourier series, that
is, series which are not necessarily absolutely and uniformly convergent,
with analytic functions, but rather with at most piecewise continuous and
differentiable functions that can be non-differentiable or even
discontinuous at some points. Therefore, it comes as a bit of a surprise
that {\em all} such series, if convergent, have limiting functions that
are also limits of inner analytic functions on the open unit disk when one
approaches the boundary of this maximum disk of convergence of the
corresponding Taylor series around the origin. This implies that, so long
as there is at most a finite number of sufficiently soft singularities on
the unit circle the series in fact converge to piecewise segments of
$C^{\infty}$ functions.

Furthermore, the relation with analytic functions in the unit disk allows
one to recover the functions that generated a Fourier series from the
coefficients of the series, even if the Fourier series itself does not
converge. In principle, this can be accomplished by taking limits of the
corresponding inner analytic functions from within the unit disk to the
unit circle in the complex plane. This helps to give to the Fourier
coefficients a definite meaning even when the corresponding series are
divergent. Usually, this meaning is thought of in terms of the idea that
the coefficients of the series still identify the original function, even
if the series does not converge. Now we are able to present a more
powerful and general way to recover the function from the coefficients in
such cases.

Note that any function $w(z)$ which is analytic on the open unit disk and
satisfies the other conditions defining it as an inner analytic function
still corresponds to a pair of FC Fourier series at the unit circle, even
if these series are divergent everywhere over the unit circle. In fact, to
accomplish this correspondence it would be enough to establish that the
open unit disk is the maximum disk of convergence of the Taylor series of
the inner analytic function $w(z)$, for example using the ratio test. This
seems to indicate that the analytic-function structure on the unit disk is
the more fundamental one. The fact that we may use the analytic structure
to recover the function $f(\theta)$ from the Fourier coefficients, even
when the Fourier series is everywhere divergent, seems to indicate the
same thing.

The analysis of convergence of a given complete Fourier series will
require, in general, the separate discussion of the convergence of the
cosine and sine sub-series, and thus will involve two pairs of DP Fourier
series and two inner analytic functions $w(z)$. Whether or not each pair
of FC Fourier series converge to the corresponding functions depends on
whether or not the Taylor series of the analytic function $w(z)$ converges
on the unit circle. This is always the case if $w(z)$ is analytic over the
whole unit circle, but otherwise its Taylor series may or may not be
convergent at some or at all points on that circle. However, when and
where $S_{z}$ does converge on the unit circle, that convergence implies
the convergence of $S_{v}$ and therefore of the two corresponding FC
Fourier series.

Since the convergence of the DP Fourier series is determined by the
convergence of the power series $S_{z}$ and $S_{v}$, it also follows from
our results here that the convergence or divergence of all these series is
ruled by the singularity structure of the inner analytic functions $w(z)$.
In the two extreme cases examined here, it is the existence of absence of
singularities within the unit disk that determines everything. In the case
of strong divergence, the existence of a single singularity within the
open unit disk suffices to determine the divergence of $S_{v}$ over the
whole unit circle, and thus the divergence of the DP Fourier series almost
everywhere over the periodic interval. In the case of strong convergence,
the absence of singularities over the whole closed unit disk suffices to
establish the strong convergence of $S_{v}$ over the whole unit circle,
and thus the strong convergence of the DP Fourier series everywhere over
the periodic interval. The remaining case is that in which there are
singularities only on the unit circle, and not within the open unit disk.
This is the case in which that disk is the maximum disk of convergence of
the series $S_{z}$. This is the more complex case, and will be examined in
the follow-up paper~\cite{FTotCPII}.

It is conceivable that the relation of Fourier theory with complex
analysis presented here can be used for other ends, possibly more general
and abstract. For example, if one is looking for a necessary and/or
sufficient condition for the convergence of Fourier series, one may
consider trying to establish such a necessary and/or sufficient condition
for the convergence of Taylor series of inner analytic functions on the
unit circle, including the cases in which the unit circle is the rim of
the maximum disk of convergence of the series. Also, as we have show here
for the simple case of the Dirac delta ``function'', it is possible that
the distributions associated to such singular objects can be discussed in
the complex plane, by means of the use of their Fourier representations.

Since complex analysis and analytic functions constitute such a powerful
tool, with so many applications in almost all areas of mathematics and
physics, it is to be hoped that other applications of the ideas explored
here will in due time present themselves. It should be noted that the
relationship with analytic functions and Taylor series constitutes a new
way to present the subject of Fourier series, that in fact might become a
rather simple and straightforward way to teach the subject.

\section{Acknowledgements}

The author would like to thank his friend and colleague Prof. Carlos
Eugênio Imbassay Carneiro, to whom he is deeply indebted for all his
interest and help, as well as his careful reading of the manuscript and
helpful criticism regarding this work.

\appendix

\section{Appendix: Riemann's Function}\label{APPriemannf}

It is interesting to point out that a certain lacunary (or lacunar)
trigonometric series proposed by Riemann, as a possible example of a
continuous function that is not differentiable anywhere, can be analyzed
in the context of the relationship presented in this paper. Riemann
proposed the sine series

\begin{displaymath}
  S_{\rm s}
  =
  \sum_{k=1}^{\infty}
  \frac{1}{k^{2}}\,
  \sin\!\left(k^{2}\theta\right),
\end{displaymath}

\noindent
from which we can construct the FC cosine series $S_{\rm c}=\bar{S}_{\rm
  s}$,

\begin{displaymath}
  S_{\rm c}
  =
  \sum_{k=1}^{\infty}
  \frac{1}{k^{2}}\,
  \cos\!\left(k^{2}\theta\right).
\end{displaymath}

\noindent
Both series pass the Weierstrass $M$-test and are therefore uniformly
convergent to continuous functions over the whole periodic interval.
Their term-by-term derivatives, however, are not convergent at all,

\noindent
\begin{eqnarray*}
  \frac{dS_{\rm c}}{d\theta}
  & = &
  -
  \sum_{k=1}^{\infty}
  \sin\!\left(k^{2}\theta\right),
  \\
  \frac{dS_{\rm s}}{d\theta}
  & = &
  \sum_{k=1}^{\infty}
  \cos\!\left(k^{2}\theta\right).
\end{eqnarray*}

\noindent
All these series can be understood as trigonometric series in which the
coefficients $a_{k}$ corresponding to values of $k$ which are not perfect
squares are all zero. Therefore, there are lots of missing terms, hence
the name ``lacunary'' series. Since the two trigonometric series $S_{\rm
  c}$ and $S_{\rm s}$ are the FC series of each other, they can be used to
construct the complex power series

\begin{displaymath}
  S_{z}
  =
  \sum_{k=1}^{\infty}
  \frac{1}{k^{2}}\,
  z^{k^{2}}.
\end{displaymath}

\noindent
This is a lacunary power series which presumably has a countable infinity
of soft singularities densely distributed over the unit circle. The ratio
test tells us that the maximum disk of convergence of this series is
indeed the open unit disk. We have therefore an inner analytic function
defined on that disk,

\noindent
\begin{eqnarray*}
  w(z)
  & = &
  \sum_{k=1}^{\infty}
  \frac{1}{k^{2}}\,
  z^{k^{2}}
  \\
  & = &
  f_{\rm c}(\rho,\theta)
  +
  \ii
  f_{\rm s}(\rho,\theta),
\end{eqnarray*}

\noindent
where

\noindent
\begin{eqnarray*}
  f_{\rm c}(1,\theta)
  & = &
  S_{\rm c},
  \\
  f_{\rm s}(1,\theta)
  & = &
  S_{\rm s}.
\end{eqnarray*}

\noindent
Since the trigonometric series are convergent on the whole unit circle,
due to Abel's theorem~\cite{AbelTheo} the $\rho\to 1$ limits of this
function from within the unit disk exists everywhere on that circle. The
derivative $w'(z)$ of this function is also analytic on the open unit
disk, and is given by

\noindent
\begin{eqnarray*}
  w'(z)
  & = &
  \frac{dw(z)}{dz}
  \\
  & = &
  \sum_{k=1}^{\infty}
  z^{k^{2}-1}.
\end{eqnarray*}

\noindent
The non-differentiability of the function $w(z)$ on the unit circle must
therefore be related to the non-existence of the $\rho\to 1$ limits of the
function $w'(z)$ from within the unit disk. Note that in terms of
derivatives with respect to $\theta$ we may write

\noindent
\begin{eqnarray*}
  \frac{dw(z)}{d\theta}
  & = &
  \ii
  z\,
  \frac{dw(z)}{dz}
  \\
  & = &
  \ii
  \sum_{k=1}^{\infty}
  z^{k^{2}}.
\end{eqnarray*}

\noindent
This relation holds true everywhere within the open unit disk. One might
therefore consider analyzing the differentiability of $f_{\rm
  c}(\rho,\theta)$ and $f_{\rm s}(\rho,\theta)$ on the unit circle by
examining the behavior of the $\rho\to 1$ limits of the function $\ii
zw'(z)$, which is proportional to the logarithmic derivative of $w(z)$.
The expectation is that the function $w(z)$ has a densely distributed
countable infinity of borderline soft singularities over the unit circle.
It then follows that its logarithmic derivative has a densely distributed
countable infinity of borderline hard singularities over the unit circle,
and therefore fails to be defined at all points. This can be understood as
a consequence of the extended version of Abel's theorem~\cite{AbelTheo},
so long as it can be shown that the series diverges to infinity everywhere
over the unit circle.

The Riemann sine series has been widely studied~\cite{jgerver}, but its FC
series less so, it would seem. There are many papers and results about
lacunary functions related to lacunary power series~\cite{LacuFunc}, which
could be of use here. In particular the Fabry and the Ostrowski-Hadamard
gap theorems~\cite{GapTheor} seem to be directly relevant. Such a lacunary
function can be characterized as one that cannot be analytically continued
beyond the boundary of the maximum disk of convergence of its Taylor
series. This would be the case if such a function had a densely
distributed set of singularities over the unit circle. Note that using the
correspondence established in this paper we may construct from any
lacunary trigonometric series a corresponding lacunary power series, and
vice-versa, so that these two subjects are completely identified with each
other.

\section{Appendix: Technical Proofs}\label{APPproofs}

\subsection{Evaluations of Convergence}\label{APPevalconv}

\begin{figure}[ht]
  \centering
  \fbox{
    \epsfig{file=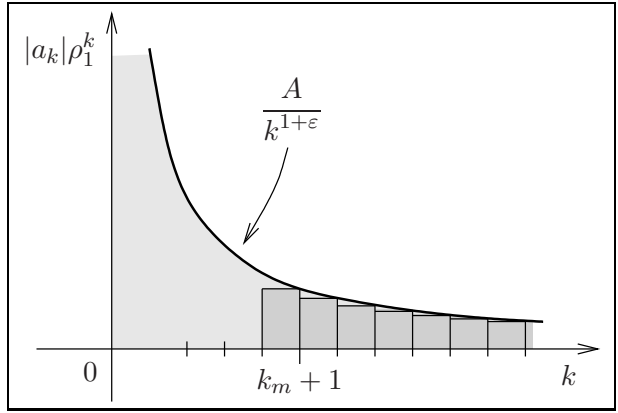,scale=1.0,angle=0}
  }
  \caption{Illustration of the upper bounding of the sum of the terms
    $|a_{k}|\rho_{1}^{k}$ by an integral.}
  \label{APPfig01}
\end{figure}

It is not a difficult task to establish the absolute and uniform
convergence of a complex power series on a circle centered at the origin,
or the lack thereof, starting from the behavior of the terms of the series
in the limit $k\to\infty$, if we assume that they behave as inverse powers
of $k$ for large values of $k$. If we have a complex power series $S_{z}$
with coefficients $a_{k}$, at a point $z_{1}$ strictly inside the unit
disk,

\begin{displaymath}
  S_{z}
  =
  \sum_{k=1}^{\infty}a_{k}z_{1}^{k},
\end{displaymath}

\noindent
where $|z_{1}|=\rho_{1}$ and $\rho_{1}<1$, then it is absolutely
convergent if and only if the series $\overline{S}_{z}$ of the absolute
values of the coefficients,

\noindent
\begin{eqnarray*}
  \overline{S}_{z}
  & = &
  \sum_{k=1}^{\infty}|a_{k}||z_{1}|^{k}
  \\
  & = &
  \sum_{k=1}^{\infty}|a_{k}|\rho_{1}^{k},
\end{eqnarray*}

\noindent
converges. One can show that this sum will be finite if, for $k$ above a
certain minimum value $k_{m}$, it holds that

\begin{displaymath}
  |a_{k}|\rho_{1}^{k}
  \leq
  \frac{A}{k^{1+\varepsilon}},
\end{displaymath}

\noindent
for some positive real constant $A$ and some real constant
$\varepsilon>0$. This is true because the sum of a finite set of initial
terms is necessarily finite, and because in this case we may bound the
remaining infinite sum from above by a convergent asymptotic integral,

\noindent
\begin{eqnarray*}
  \sum_{k=k_{m}+1}^{\infty}|a_{k}|\rho_{1}^{k}
  & \leq &
  \sum_{k=k_{m}+1}^{\infty}\frac{A}{k^{1+\varepsilon}}
  \\
  & < &
  \int_{k_{m}}^{\infty}dk\,\frac{A}{k^{1+\varepsilon}}
  \\
  & = &
  \frac{-A}{\varepsilon}\,\frac{1}{k^{\varepsilon}}\at{k_{m}}{\infty}
  \\
  & = &
  \frac{A}{\varepsilon}\,\frac{1}{k_{m}^{\varepsilon}},
\end{eqnarray*}

\begin{figure}[ht]
  \centering
  \fbox{
    \epsfig{file=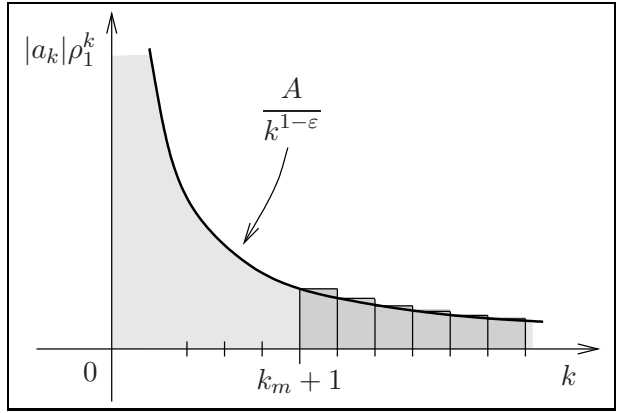,scale=1.0,angle=0}
  }
  \caption{Illustration of the lower bounding of the sum of the terms
    $|a_{k}|\rho_{1}^{k}$ by an integral.}
  \label{APPfig02}
\end{figure}

\noindent
as illustrated in Figure~\ref{APPfig01}. In that illustration each
vertical rectangle has base $1$ and height given by $|a_{k}|\rho_{1}^{k}$,
and therefore area given by $|a_{k}|\rho_{1}^{k}$. As one can see, the
construction is such that the set of all such rectangles is below the
graph of the function $A/k^{1+\varepsilon}$, and therefore the sum of
their areas is contained within the area under that graph, to the right of
$k_{m}$. This establishes the necessary inequality between the sum and the
integral.

So long as $\varepsilon$ is not zero, this establishes an upper bound to a
sum of positive quantities, which is therefore a monotonically increasing
sum. It then follows from the well-known theorem of real analysis that the
sum necessarily converges, and therefore the series $S_{z}$ is absolutely
convergent at $z_{1}$.

In addition to this, one can see that the convergence condition does not
depend on $\theta$, since that dependence is only within the complex
variable $z=\rho_{1}\exp(\ii\theta)$, and vanishes when we take absolute
values. This implies uniform convergence because, given a strictly
positive real number $\epsilon$, absolute convergence for this value of
$\epsilon$ implies convergence for this same value of $\epsilon$, with the
same solution $k(\epsilon)$ for the convergence condition. This makes it
clear that the solution of the convergence condition for $k$ is
independent of position and therefore that the series is also uniformly
convergent on the circle of radius $\rho_{1}$ centered at the origin.

This establishes a sufficient condition for the absolute and uniform
convergence of the complex power series on the circle. On the other hand,
if we have that, for $k$ above a certain minimum value $k_{m}$,

\begin{displaymath}
  |a_{k}|\rho_{1}^{k}
  \geq
  \frac{A}{k^{1-\varepsilon}},
\end{displaymath}

\noindent
with positive real $A$ and real $\varepsilon\geq 0$, then it is possible
to bound the sum $\overline{S}_{v}$ from below by an asymptotic integral
that diverges to positive infinity. This is done in a way similar to the
one used for the establishment of the upper bound, but inverting the
situation so as to keep the area under the graph contained within the
combined areas of the rectangles, as illustrated in Figure~\ref{APPfig02}.
The argument then establishes in this case that, for $\varepsilon>0$

\noindent
\begin{eqnarray*}
  \sum_{k=k_{m}+1}^{\infty}|a_{k}|\rho_{1}^{k}
  & \geq &
  \sum_{k=k_{m}+1}^{\infty}\frac{A}{k^{1-\varepsilon}}
  \\
  & > &
  \int_{k_{m}+1}^{\infty}dk\,\frac{A}{k^{1-\varepsilon}}
  \\
  & = &
  \frac{A}{\varepsilon}\,k^{\varepsilon}\at{k_{m}+1}{\infty}
  \\
  & = &
  -\,
  \frac{A}{\varepsilon}\,(k_{m}+1)^{\varepsilon}
  +
  \frac{A}{\varepsilon}
  \lim_{k\to\infty}
  k^{\varepsilon},
\end{eqnarray*}

\noindent
and therefore that $\overline{S}_{z}$ diverges to infinity. A similar
calculation can be performed in the case $\varepsilon=0$, leading to
logarithms and yielding the same conclusions. This does not prove or
disprove convergence itself, but it does establish the absence of absolute
convergence. It also shows that, so long as $|a_{k}|\rho_{1}^{k}$ behaves
as a power of $k$ for large $k$, the previous condition is both sufficient
and necessary for absolute convergence.

\subsection{About the Ratio Test}\label{APPratitest}

Let us show that the Taylor series of an analytic function cannot satisfy
the ratio test at the boundary of its maximum convergence disk. We will
assume that the convergence disk is the open unit disk, but the proof can
be easily generalized. It the series satisfies the test, then it follows
that there is a strictly positive real number $q<1$ and an integer $k_{m}$
such that, for all $k>k_{m}$ we have

\begin{displaymath}
  \frac{|a_{k+1}|}{|a_{k}|}\,|z|
  <
  q.
\end{displaymath}

\noindent
We observe now that since $|z|=1$ over the whole unit circle, this
condition is in fact independent of the angular position and valid over
the whole unit circle, so that we have

\begin{displaymath}
  \frac{|a_{k+1}|}{|a_{k}|}
  <
  q.
\end{displaymath}

\noindent
Now, consider a point $z_{1}$ given by any value of $\theta$ and the value

\begin{displaymath}
  |z_{1}|
  =
  \frac{q+1}{2q}
\end{displaymath}

\noindent
for the radius. Since $q$ is strictly positive and strictly smaller than
one, we have that

\noindent
\begin{eqnarray*}
  q
  & < &
  1
  \;\;\;\Rightarrow
  \\
  q+q
  & < &
  q+1
  \;\;\;\Rightarrow
  \\
  \frac{q+1}{2q}
  & > &
  1.
\end{eqnarray*}

\noindent
It follows therefore that the point $z_{1}$ is strictly outside the unit
disk. We consider now the ratio test for this new point. We start from the
known valid condition on the unit disk, so that we have

\noindent
\begin{eqnarray*}
  \frac{|a_{k+1}|}{|a_{k}|}
  & < &
  q
  \;\;\;\Rightarrow
  \\
  \frac{|a_{k+1}|}{|a_{k}|}\,|z_{1}|
  & < &
  \frac{q+1}{2}
  \;\;\;\Rightarrow
  \\
  \frac{|a_{k+1}|}{|a_{k}|}\,|z_{1}|
  & < &
  q_{1},
\end{eqnarray*}

\noindent
which establishes the upper bound $q_{1}$. Besides, since we know that $q$
is strictly less than one, we also have that

\noindent
\begin{eqnarray*}
  q
  & < &
  1
  \;\;\;\Rightarrow
  \\
  q+1
  & < &
  1+1
  \;\;\;\Rightarrow
  \\
  \frac{q+1}{2}
  & < &
  1.
\end{eqnarray*}

\noindent
What we have concluded here is that $q_{1}<1$, and therefore we conclude
that the series satisfies the ratio test at the new point, which is
strictly outside the unit disk, and hence converges there. Therefore, the
basic convergence theorem implies that the maximum disk of convergence of
the series $S_{z}$ extends beyond the unit circle. This contradicts the
hypothesis that the open unit disk is the maximum disk of convergence, and
we conclude therefore that the series cannot satisfy the ratio test at the
boundary of its maximum disk of convergence.

Another way to interpret this result is to say that, if the Taylor series
$S_{z}$ does satisfy the ratio test at the boundary of a given disk, then
that disk is not its maximum disk of convergence. It follows therefore
that the maximum disk of convergence is larger than the given disk, and
contains it.

\section{Appendix: Examples of Limits from Within}\label{APPexamplim}

In this appendix we will give a few illustrative examples. We will explore
the possibility of recovering values for the functions associated to
either weakly convergent or flatly divergent Fourier series, using limits
from the interior of the open unit disk. In some of these simple examples
we will recover the original functions in closed form. As we will see,
sometimes it is also possible to write the FC functions in closed form. At
the same time, we will also explore the extension of the structure to
singular objects such as the Dirac delta ``function''.

\subsection{A Regular Sine Series with All {\boldmath $k$}}

Consider the Fourier series of the one-cycle unit-amplitude sawtooth wave,
which is just the linear function $\theta/\pi$ between $-\pi$ and $\pi$.
As is well known it is given by the sine series

\begin{displaymath}
  S_{\rm s}
  =
  -\,
  \frac{2}{\pi}
  \sum_{k=1}^{\infty}
  \frac{(-1)^{k}}{k}\,
  \sin(k\theta).
\end{displaymath}

\noindent
The corresponding FC series is then

\begin{displaymath}
  \bar{S}_{\rm s}
  =
  -\,
  \frac{2}{\pi}
  \sum_{k=1}^{\infty}
  \frac{(-1)^{k}}{k}\,
  \cos(k\theta),
\end{displaymath}

\noindent
the complex $S_{v}$ series is given by

\begin{displaymath}
  S_{v}
  =
  -\,
  \frac{2}{\pi}
  \sum_{k=1}^{\infty}
  \frac{(-1)^{k}}{k}\,
  v^{k},
\end{displaymath}

\noindent
and the complex power series $S_{z}$ is given by

\begin{displaymath}
  S_{z}
  =
  -\,
  \frac{2}{\pi}
  \sum_{k=1}^{\infty}
  \frac{(-1)^{k}}{k}\,
  z^{k}.
\end{displaymath}

\noindent
The ratio test tells us that the disk of convergence of $S_{z}$ is the
unit disk. If we consider the inner analytic function $w(z)$ within this
disk we observe that $w(0)=0$, as expected. We have for the function and
its derivative

\noindent
\begin{eqnarray*}
  w(z)
  & = &
  -\,
  \frac{2}{\pi}
  \sum_{k=1}^{\infty}
  \frac{(-1)^{k}}{k}\,
  z^{k}
  \;\;\;\Rightarrow
  \\
  \frac{dw(z)}{dz}
  & = &
  -\,
  \frac{2}{\pi}
  \sum_{k=1}^{\infty}
  (-1)^{k}
  z^{k-1}
  \\
  & = &
  -\,
  \frac{2}{\pi z}
  \sum_{k=1}^{\infty}
  (-z)^{k}.
\end{eqnarray*}

\noindent
This is a geometrical series and therefore we may write the derivative in
closed form,

\noindent
\begin{eqnarray*}
  \frac{dw(z)}{dz}
  & = &
  -\,
  \frac{2}{\pi z}\,
  \frac{-z}{1+z}
  \\
  & = &
  \frac{2}{\pi}\,
  \frac{1}{1+z}.
\end{eqnarray*}

\noindent
This function has a simple pole at the point where the one-cycle sawtooth
wave is discontinuous. This is a hard singularity of degree $1$. It is now
very simple to integrate in order to obtain the inner analytic function
$w(z)$ in closed form, remembering that we should have $w(0)=0$,

\begin{displaymath}
  w(z)
  =
  \frac{2}{\pi}\,
  \ln(1+z).
\end{displaymath}

\noindent
This function has a logarithmic singularity at the point where the
one-cycle sawtooth wave is discontinuous. This is a borderline hard
singularity. If we now rationalize the argument of the logarithm in order
to write it in the form

\begin{displaymath}
  1+z
  =
  R\e{\ii\alpha},
\end{displaymath}

\noindent
then the logarithm is given by

\begin{displaymath}
  \ln\!\left(1+z\right)
  =
  \ln(R)+\ii\alpha,
\end{displaymath}

\noindent
which allows us to identify the real and imaginary parts of $w(z)$. We may
write for the argument

\noindent
\begin{eqnarray*}
  R\e{\ii\alpha}
  & = &
  1
  +
  \rho
  \left[
    \cos(\theta)
    +
    \ii
    \sin(\theta)
  \right]
  \\
  & = &
  \left[
    1
    +
    \rho
    \cos(\theta)
  \right]
  +
  \ii
  \left[
    \rho
    \sin(\theta)
  \right].
\end{eqnarray*}

\noindent
From this it follows, after some algebra, that we have

\noindent
\begin{eqnarray*}
  R(\rho,\theta)
  & = &
  \sqrt{1+2\rho\cos(\theta)+\rho^{2}},
  \\
  \alpha(\rho,\theta)
  & = &
  \arctan\!
  \left[\frac{\rho\sin(\theta)}{1+\rho\cos(\theta)}\right].
\end{eqnarray*}

\noindent
The part which is of interest now is the imaginary one, which is related
to the series $S_{\rm s}$,

\begin{displaymath}
  \Im[w(z)]
  =
  \frac{2}{\pi}\,
  \arctan\!
  \left[\frac{\rho\sin(\theta)}{1+\rho\cos(\theta)}\right].
\end{displaymath}

\noindent
The $\rho\to 1$ limit of this quantity is the function $f_{\rm
  s}(\theta)$. If we consider the case $\theta=\pm\pi$, we immediately get
zero, because the argument of the arc tangent is zero, and therefore we
get $\alpha=0$. Therefore we have

\noindent
\begin{eqnarray*}
  f_{\rm s}(\pm\pi)
  & = &
  \lim_{\rho\to 1}
  \Im[w(\rho,\pm\pi)]
  \\
  & = &
  0.
\end{eqnarray*}

\noindent
This is the correct value for the Fourier-series representation of the
one-cycle sawtooth wave at this point. For $\theta\neq\pm\pi$ we get

\noindent
\begin{eqnarray*}
  f_{\rm s}(\theta)
  & = &
  \lim_{\rho\to 1}
  \Im[w(z)]
  \\
  & = &
  \frac{2}{\pi}
  \lim_{\rho\to 1}
  \arctan\!
  \left[\frac{\rho\sin(\theta)}{1+\rho\cos(\theta)}\right]
  \\
  & = &
  \frac{2}{\pi}
  \arctan\!
  \left[\frac{\sin(\theta)}{1+\cos(\theta)}\right].
\end{eqnarray*}

\noindent
In order to solve this we write everything in terms of $\theta/2$, and
thus obtain

\noindent
\begin{eqnarray*}
  f_{\rm s}(\theta)
  & = &
  \frac{2}{\pi}
  \arctan\!
  \left[\frac{2\sin(\theta/2)\cos(\theta/2)}{2\cos^{2}(\theta/2)}\right]
  \\
  & = &
  \frac{2}{\pi}
  \arctan\!\left[\frac{\sin(\theta/2)}{\cos(\theta/2)}\right]
  \\
  & = &
  \frac{2}{\pi}
  \arctan\!\left[\tan(\theta/2)\right]
  \\
  & = &
  \frac{\theta}{\pi}.
\end{eqnarray*}

\noindent
This is the correct value for the one-cycle unit-amplitude sawtooth wave.
The FC function is related to the real part

\begin{displaymath}
  \Re[w(z)]
  =
  \frac{1}{\pi}\,
  \ln\!\left[1+2\rho\cos(\theta)+\rho^{2}\right].
\end{displaymath}

\noindent
In the $\rho\to 1$ limit this function is singular at $\theta=\pm\pi$.
Away from this logarithmic singularity we may write

\noindent
\begin{eqnarray*}
  \bar{f}_{\rm s}(\theta)
  & = &
  \lim_{\rho\to 1}
  \Re[w(z)]
  \\
  & = &
  \frac{1}{\pi}\,
  \ln\!\left[2+2\cos(\theta)\right]
  \\
  & = &
  \frac{1}{\pi}\,
  \ln\!\left[4\cos^{2}(\theta/2)\right].
\end{eqnarray*}

\noindent
We see therefore that the FC function to the one-cycle sawtooth wave is a
function with a logarithmic singularity at the point where the original
function is discontinuous. The derivative of this conjugate function is
easily calculated, and turns out to be

\begin{displaymath}
  \frac{d\bar{f}_{\rm s}(\theta)}{d\theta}
  =
  \frac{1}{\pi}\,
  \frac{\cos(\theta)-1}{\sin(\theta)}.
\end{displaymath}

\noindent
This is a quite a simple function, with a simple pole at $\theta=\pm\pi$.
Note that the logarithmic singularity of the FC function $\bar{f}_{\rm
  s}(\theta)$ is an integrable one, as one would expect since the Fourier
coefficients are finite.

Note that we get a free set of integration formulas out of this effort.
Since the coefficients $a_{k}$ can be written as integrals involving
$\bar{f}_{\rm s}(\theta)$, we have at once that

\begin{displaymath}
  \int_{-\pi}^{\pi}d\theta\,
  \ln\!\left[4\cos^{2}(\theta/2)\right]
  \cos(k\theta)
  =
  -\,
  \frac{2\pi(-1)^{k}}{k},
\end{displaymath}

\noindent
for $k>0$.

\subsection{A Regular Sine Series with Odd {\boldmath
    $k$}}\label{APPsquarewave}

Consider the Fourier series of the standard unit-amplitude square wave. As
is well known it is given by the sine series

\begin{displaymath}
  S_{\rm s}
  =
  \frac{4}{\pi}
  \sum_{j=0}^{\infty}
  \frac{1}{2j+1}\,
  \sin[(2j+1)\theta].
\end{displaymath}

\noindent
The corresponding FC series is then

\begin{displaymath}
  \bar{S}_{\rm s}
  =
  \frac{4}{\pi}
  \sum_{j=0}^{\infty}
  \frac{1}{2j+1}\,
  \cos[(2j+1)\theta],
\end{displaymath}

\noindent
the complex $S_{v}$ series is given by

\begin{displaymath}
  S_{v}
  =
  \frac{4}{\pi}
  \sum_{j=0}^{\infty}
  \frac{1}{2j+1}\,
  v^{2j+1},
\end{displaymath}

\noindent
and the complex power series $S_{z}$ is given by

\begin{displaymath}
  S_{z}
  =
  \frac{4}{\pi}
  \sum_{j=0}^{\infty}
  \frac{1}{2j+1}\,
  z^{2j+1}.
\end{displaymath}

\noindent
The ratio test tells us that the disk of convergence of $S_{z}$ is the
unit disk. If we consider the inner analytic function $w(z)$ within this
disk we observe that $w(0)=0$, as expected. We have for the function and
its derivative

\noindent
\begin{eqnarray*}
  w(z)
  & = &
  \frac{4}{\pi}
  \sum_{j=0}^{\infty}
  \frac{1}{2j+1}\,
  z^{2j+1}
  \;\;\;\Rightarrow
  \\
  \frac{dw(z)}{dz}
  & = &
  \frac{4}{\pi}
  \sum_{j=0}^{\infty}
  z^{2j}.
\end{eqnarray*}

\noindent
This is a geometrical series and therefore we may write the derivative in
closed form,

\noindent
\begin{eqnarray*}
  \frac{dw(z)}{dz}
  & = &
  \frac{4}{\pi}\,
  \frac{1}{1-z^{2}}
  \\
  & = &
  \frac{2}{\pi}\,
  \left(
    \frac{1}{1+z}
    +
    \frac{1}{1-z}
  \right).
\end{eqnarray*}

\noindent
This function has two simple poles at the points where the square wave is
discontinuous. These are hard singularities of degree $1$. It is now very
simple to integrate in order to obtain the inner analytic function $w(z)$
in closed form, remembering that we should have $w(0)=0$,

\noindent
\begin{eqnarray*}
  w(z)
  & = &
  \frac{2}{\pi}\,
  \left[
    \ln(1+z)
    -
    \ln(1-z)
  \right]
  \\
  & = &
  \frac{2}{\pi}\,
  \ln\!\left(\frac{1+z}{1-z}\right).
\end{eqnarray*}

\noindent
This function has logarithmic singularities at the points where the square
wave is discontinuous. These are borderline hard singularities. If we now
rationalize the argument of the logarithm in order to write it in the form

\begin{displaymath}
  \frac{1+z}{1-z}
  =
  R\e{\ii\alpha},
\end{displaymath}

\noindent
then the logarithm is given by

\begin{displaymath}
  \ln\!\left(\frac{1+z}{1-z}\right)
  =
  \ln(R)+\ii\alpha,
\end{displaymath}

\noindent
which allows us to identify the real and imaginary parts of $w(z)$. We may
write for the argument

\noindent
\begin{eqnarray*}
  R\e{\ii\alpha}
  & = &
  \frac{(1+z)(1-z*)}{(1-z)(1-z*)}
  \\
  & = &
  \frac
  {1+2\ii\rho\sin(\theta)-\rho^{2}}
  {1-2\rho\cos(\theta)+\rho^{2}}
  \\
  & = &
  \frac
  {\left(1-\rho^{2}\right)}
  {\left(1+\rho^{2}\right)-2\rho\cos(\theta)}
  +
  \ii\,
  \frac
  {2\rho\sin(\theta)}
  {\left(1+\rho^{2}\right)-2\rho\cos(\theta)}.
\end{eqnarray*}

\noindent
From this it follows, after some algebra, that we have

\noindent
\begin{eqnarray*}
  R(\rho,\theta)
  & = &
  \sqrt
  {
    \frac
    {\left(1+\rho^{2}\right)+2\rho\cos(\theta)}
    {\left(1+\rho^{2}\right)-2\rho\cos(\theta)}
  },
  \\
  \alpha(\rho,\theta)
  & = &
  \arctan\!\left[\frac{2\rho\sin(\theta)}{1-\rho^{2}}\right].
\end{eqnarray*}

\noindent
The part which is of interest now is the imaginary one, which is related
to the series $S_{\rm s}$,

\begin{displaymath}
  \Im[w(z)]
  =
  \frac{2}{\pi}\,
  \arctan\!\left[\frac{2\rho\sin(\theta)}{1-\rho^{2}}\right].
\end{displaymath}

\noindent
The $\rho\to 1$ limit of this quantity is the function $f_{\rm
  s}(\theta)$. If we consider the cases $\theta=0$ and $\theta=\pm\pi$, we
immediately get zero, because the argument of the arc tangent is zero, and
therefore we get $\alpha=0$. Therefore we have

\noindent
\begin{eqnarray*}
  f_{\rm s}(0)
  & = &
  \lim_{\rho\to 1}
  \Im[w(\rho,0)]
  \\
  & = &
  0,
  \\
  f_{\rm s}(\pm\pi)
  & = &
  \lim_{\rho\to 1}
  \Im[w(\rho,\pm\pi)]
  \\
  & = &
  0.
\end{eqnarray*}

\noindent
These are the correct values for the Fourier-series representation of the
square wave at these points. For $\theta\neq 0$ and $\theta\neq\pm\pi$ we
get

\noindent
\begin{eqnarray*}
  f_{\rm s}(\theta)
  & = &
  \lim_{\rho\to 1}
  \Im[w(z)]
  \\
  & = &
  \frac{2}{\pi}
  \lim_{\rho\to 1}
  \arctan\!\left[\frac{2\rho\sin(\theta)}{1-\rho^{2}}\right]
  \\
  & = &
  \frac{2}{\pi}
  \lim_{\rho\to 1}
  \arctan\!\left[\frac{2\sin(\theta)}{1-\rho^{2}}\right].
\end{eqnarray*}

\noindent
If $\theta>0$ the argument of the arc tangent goes to positive infinity in
the limit, and therefore the arc tangent $\alpha$ approaches $\pi/2$. If
$\theta>0$ the argument goes to negative infinity, and therefore $\alpha$
approaches $-\pi/2$. We get therefore the values

\noindent
\begin{displaymath}
  %
  \renewcommand{\arraystretch}{2.0}
  \begin{array}{rclcc}
    f_{\rm s}(\theta)
    & = &
    1
    &
    \mbox{for}
    &
    \theta>0,
    \\
    f_{\rm s}(\theta)
    & = &
    -1
    &
    \mbox{for}
    &
    \theta<0,
  \end{array}
\end{displaymath}

\noindent
which completes the correct set of values for the standard square wave.
The FC function is related to the real part

\begin{displaymath}
  \Re[w(z)]
  =
  \frac{1}{\pi}\,
  \ln\!
  \left[
    \frac
    {\left(1+\rho^{2}\right)+2\rho\cos(\theta)}
    {\left(1+\rho^{2}\right)-2\rho\cos(\theta)}
  \right].
\end{displaymath}

\noindent
In the $\rho\to 1$ limit this function is singular at $\theta=0$ and
$\theta=\pm\pi$. Away from these logarithmic singularities we may write

\noindent
\begin{eqnarray*}
  \bar{f}_{\rm s}(\theta)
  & = &
  \lim_{\rho\to 1}
  \Re[w(z)]
  \\
  & = &
  \frac{1}{\pi}\,
  \ln\!\left[\frac{1+\cos(\theta)}{1-\cos(\theta)}\right].
\end{eqnarray*}

\noindent
We see therefore that the FC function to the square wave is a function
with logarithmic singularities at the two points where the original
function is discontinuous. The derivative of this conjugate function is
easily calculated, and turns out to be

\begin{displaymath}
  \frac{d\bar{f}_{\rm s}(\theta)}{d\theta}
  =
  -\,
  \frac{2}{\pi}\,
  \frac{1}{\sin(\theta)}.
\end{displaymath}

\noindent
This is a quite a simple function, with two simple poles at $\theta=0$ and
$\theta=\pm\pi$. Note that the logarithmic singularities of the FC
function $\bar{f}_{\rm s}(\theta)$ are integrable ones, as one would
expect since the Fourier coefficients are finite.

Note that we get a free set of integration formulas out of this effort.
Since the coefficients $a_{k}$ can be written as integrals involving
$\bar{f}_{\rm s}(\theta)$, we have at once that for odd $k=2j+1$

\begin{displaymath}
  \int_{-\pi}^{\pi}d\theta\,
  \ln\!\left[\frac{1+\cos(\theta)}{1-\cos(\theta)}\right]
  \cos(k\theta)
  =
  \frac{4\pi}{k},
\end{displaymath}

\noindent
while for even $k=2j$

\begin{displaymath}
  \int_{-\pi}^{\pi}d\theta\,
  \ln\!\left[\frac{1+\cos(\theta)}{1-\cos(\theta)}\right]
  \cos(k\theta)
  =
  0,
\end{displaymath}

\noindent
for $k>0$. Note that $\bar{f}_{\rm s}(\theta)$ is an even function of
$\theta$, so that these integrals are not zero by parity arguments, and
therefore we also have

\begin{displaymath}
  \int_{0}^{\pi}d\theta\,
  \ln\!\left[\frac{1+\cos(\theta)}{1-\cos(\theta)}\right]
  \cos(k\theta)
  =
  0,
\end{displaymath}

\noindent
for $k>0$.

\subsection{A Regular Sine Series with Even {\boldmath $k$}}

Consider the Fourier series of the two-cycle unit-amplitude sawtooth wave.
As is well known it is given by the sine series

\begin{displaymath}
  S_{\rm s}
  =
  -\,
  \frac{4}{\pi}
  \sum_{j=1}^{\infty}
  \frac{1}{2j}\,
  \sin[(2j)\theta].
\end{displaymath}

\noindent
The corresponding FC series is then

\begin{displaymath}
  \bar{S}_{\rm s}
  =
  -\,
  \frac{4}{\pi}
  \sum_{j=1}^{\infty}
  \frac{1}{2j}\,
  \cos[(2j)\theta],
\end{displaymath}

\noindent
the complex $S_{v}$ series is given by

\begin{displaymath}
  S_{v}
  =
  -\,
  \frac{4}{\pi}
  \sum_{j=1}^{\infty}
  \frac{1}{2j}\,
  v^{2j},
\end{displaymath}

\noindent
and the complex power series $S_{z}$ is given by

\begin{displaymath}
  S_{z}
  =
  -\,
  \frac{4}{\pi}
  \sum_{j=1}^{\infty}
  \frac{1}{2j}\,
  z^{2j}.
\end{displaymath}

\noindent
The ratio test tells us that the disk of convergence of $S_{z}$ is the
unit disk. If we consider the inner analytic function $w(z)$ within this
disk we observe that $w(0)=0$, as expected. We have for the function and
its derivative

\noindent
\begin{eqnarray*}
  w(z)
  & = &
  -\,
  \frac{4}{\pi}
  \sum_{j=1}^{\infty}
  \frac{1}{2j}\,
  z^{2j}
  \;\;\;\Rightarrow
  \\
  \frac{dw(z)}{dz}
  & = &
  -\,
  \frac{4}{\pi}
  \sum_{j=1}^{\infty}
  z^{2j-1}.
\end{eqnarray*}

\noindent
This is a geometrical series and therefore we may write the derivative in
closed form,

\noindent
\begin{eqnarray*}
  \frac{dw(z)}{dz}
  & = &
  -\,
  \frac{4}{\pi}\,
  \frac{z}{1-z^{2}}
  \\
  & = &
  \frac{2}{\pi}\,
  \left(
    \frac{1}{1+z}
    -
    \frac{1}{1-z}
  \right).
\end{eqnarray*}

\noindent
This function has two simple poles at the points where the two-cycle
sawtooth wave is discontinuous. These are hard singularities of degree
$1$. It is now very simple to integrate in order to obtain the inner
analytic function $w(z)$ in closed form, remembering that we should have
$w(0)=0$,

\noindent
\begin{eqnarray*}
  w(z)
  & = &
  \frac{2}{\pi}\,
  \left[
    \ln(1+z)
    +
    \ln(1-z)
  \right]
  \\
  & = &
  \frac{2}{\pi}\,
  \ln\!\left(1-z^{2}\right).
\end{eqnarray*}

\noindent
This function has logarithmic singularities at the points where the
two-cycle sawtooth wave is discontinuous. These are borderline hard
singularities. If we now rationalize the argument of the logarithm in
order to write it in the form

\begin{displaymath}
  1-z^{2}
  =
  R\e{\ii\alpha},
\end{displaymath}

\noindent
then the logarithm is given by

\begin{displaymath}
  \ln\!\left(1-z^{2}\right)
  =
  \ln(R)+\ii\alpha,
\end{displaymath}

\noindent
which allows us to identify the real and imaginary parts of $w(z)$. We may
write for the argument

\noindent
\begin{eqnarray*}
  R\e{\ii\alpha}
  & = &
  1
  -
  \rho^{2}
  \left[
    \cos(2\theta)
    +
    \ii
    \sin(2\theta)
  \right]
  \\
  & = &
  \left[
    1
    -
    \rho^{2}
    \cos(2\theta)
  \right]
  -
  \ii
  \left[
    \rho^{2}
    \sin(2\theta)
  \right].
\end{eqnarray*}

\noindent
From this it follows, after some algebra, that we have

\noindent
\begin{eqnarray*}
  R(\rho,\theta)
  & = &
  \sqrt{1-2\rho^{2}\cos(2\theta)+\rho^{4}},
  \\
  \alpha(\rho,\theta)
  & = &
  -\arctan\!
  \left[\frac{\rho^{2}\sin(2\theta)}{1-\rho^{2}\cos(2\theta)}\right].
\end{eqnarray*}

\noindent
The part which is of interest now is the imaginary one, which is related
to the series $S_{\rm s}$,

\begin{displaymath}
  \Im[w(z)]
  =
  -\,
  \frac{2}{\pi}\,
  \arctan\!
  \left[\frac{\rho^{2}\sin(2\theta)}{1-\rho^{2}\cos(2\theta)}\right].
\end{displaymath}

\noindent
The $\rho\to 1$ limit of this quantity is the function $f_{\rm
  s}(\theta)$. If we consider the cases $\theta=0$ and $\theta=\pm\pi$, we
immediately get zero, because the argument of the arc tangent is zero, and
therefore we get $\alpha=0$. Therefore we have

\noindent
\begin{eqnarray*}
  f_{\rm s}(0)
  & = &
  \lim_{\rho\to 1}
  \Im[w(\rho,0)]
  \\
  & = &
  0,
  \\
  f_{\rm s}(\pm\pi)
  & = &
  \lim_{\rho\to 1}
  \Im[w(\rho,\pm\pi)]
  \\
  & = &
  0.
\end{eqnarray*}

\noindent
These are the correct values for the Fourier-series representation of the
two-cycle sawtooth wave at these points. For $\theta\neq 0$ and
$\theta\neq\pm\pi$ we get

\noindent
\begin{eqnarray*}
  f_{\rm s}(\theta)
  & = &
  \lim_{\rho\to 1}
  \Im[w(z)]
  \\
  & = &
  -\,
  \frac{2}{\pi}
  \lim_{\rho\to 1}
  \arctan\!
  \left[\frac{\rho^{2}\sin(2\theta)}{1-\rho^{2}\cos(2\theta)}\right]
  \\
  & = &
  -\,
  \frac{2}{\pi}
  \arctan\!
  \left[\frac{\sin(2\theta)}{1-\cos(2\theta)}\right]
  \\
  & = &
  -\,
  \frac{2}{\pi}
  \arctan\!
  \left[\frac{2\sin(\theta)\cos(\theta)}{2\sin^{2}(\theta)}\right]
  \\
  & = &
  -\,
  \frac{2}{\pi}
  \arctan\!\left[\frac{\cos(\theta)}{\sin(\theta)}\right].
\end{eqnarray*}

\noindent
In order to solve this equation we first write it as

\noindent
\begin{eqnarray*}
  f_{\rm s}(\theta)
  & = &
  \frac{2}{\pi}\,
  \alpha(\theta),
  \\
  \alpha(\theta)
  & = &
  -\,
  \arctan\!\left[\frac{\cos(\theta)}{\sin(\theta)}\right]
  \;\;\;\Rightarrow
  \\
  \tan(\alpha)
  & = &
  -\,
  \frac{1}{\tan(\theta)}.
\end{eqnarray*}

\noindent
Keeping in mind that we must have $-\pi\leq\theta\leq\pi$ and
$-\pi/2\leq\alpha\leq\pi/2$, and paying careful attention to the signs of
the sine and cosine of either angle, one concludes that the solution is
given by

\noindent
\begin{displaymath}
  %
  \renewcommand{\arraystretch}{2.4}
  \begin{array}{rclcc}
    \alpha(\theta)
    & = &
    \theta
    -
    \FFrac{\pi}{2}
    &
    \mbox{for}
    &
    \theta>0,
    \\
    \alpha(\theta)
    & = &
    \theta
    +
    \FFrac{\pi}{2}
    &
    \mbox{for}
    &
    \theta<0,
  \end{array}
\end{displaymath}

\noindent
which means that for $f_{\rm s}(\theta)$ we have

\noindent
\begin{displaymath}
  %
  \renewcommand{\arraystretch}{2.4}
  \begin{array}{rclcc}
    f_{\rm s}(\theta)
    & = &
    \FFrac{2\theta}{\pi}
    -
    1
    &
    \mbox{for}
    &
    \theta>0,
    \\
    f_{\rm s}(\theta)
    & = &
    \FFrac{2\theta}{\pi}
    +
    1
    &
    \mbox{for}
    &
    \theta<0.
  \end{array}
\end{displaymath}

\noindent
These are the correct values for the two-cycle unit-amplitude sawtooth
wave. The FC function is related to the real part

\begin{displaymath}
  \Re[w(z)]
  =
  \frac{1}{\pi}\,
  \ln\!\left[1-2\rho^{2}\cos(2\theta)+\rho^{4}\right].
\end{displaymath}

\noindent
In the $\rho\to 1$ limit this function is singular at $\theta=0$ and
$\theta=\pm\pi$. Away from these logarithmic singularities we may write

\noindent
\begin{eqnarray*}
  \bar{f}_{\rm s}(\theta)
  & = &
  \lim_{\rho\to 1}
  \Re[w(z)]
  \\
  & = &
  \frac{1}{\pi}\,
  \ln\!\left[2-2\cos(2\theta)\right]
  \\
  & = &
  \frac{1}{\pi}\,
  \ln\!\left[4\sin^{2}(\theta)\right].
\end{eqnarray*}

\noindent
We see therefore that the FC function to the two-cycle sawtooth wave is a
function with logarithmic singularities at the two points where the
original function is discontinuous. The derivative of this conjugate
function is easily calculated, and turns out to be

\begin{displaymath}
  \frac{d\bar{f}_{\rm s}(\theta)}{d\theta}
  =
  \frac{2}{\pi}\,
  \frac{\cos(\theta)}{\sin(\theta)}.
\end{displaymath}

\noindent
This is a quite a simple function, with two simple poles at $\theta=0$ and
$\theta=\pm\pi$. Note that the logarithmic singularities of the FC
function $\bar{f}_{\rm s}(\theta)$ are integrable ones, as one would
expect since the Fourier coefficients are finite.

Note that we get a free set of integration formulas out of this effort.
Since the coefficients $a_{k}$ can be written as integrals involving
$\bar{f}_{\rm s}(\theta)$, we have at once that for even $k=2j$

\begin{displaymath}
  \int_{-\pi}^{\pi}d\theta\,
  \ln\!\left[4\sin^{2}(\theta)\right]
  \cos(k\theta)
  =
  -\,
  \frac{4\pi}{k},
\end{displaymath}

\noindent
for $k>0$, while for odd $k=2j+1$

\begin{displaymath}
  \int_{-\pi}^{\pi}d\theta\,
  \ln\!\left[4\sin^{2}(\theta)\right]
  \cos(k\theta)
  =
  0.
\end{displaymath}

\noindent
Note that $\bar{f}_{\rm s}(\theta)$ is an even function of $\theta$, so
that these integrals are not zero by parity arguments, and therefore we
also have

\begin{displaymath}
  \int_{0}^{\pi}d\theta\,
  \ln\!\left[4\sin^{2}(\theta)\right]
  \cos(k\theta)
  =
  0,
\end{displaymath}

\noindent
for $k=2j+1$.

\subsection{A Regular Cosine Series with Odd {\boldmath $k$}}

Consider the Fourier series of the unit-amplitude triangular wave. As is
well known it is given by the cosine series

\begin{displaymath}
  S_{\rm c}
  =
  -\,
  \frac{8}{\pi^{2}}
  \sum_{j=0}^{\infty}
  \frac{1}{(2j+1)^{2}}\,
  \cos[(2j+1)\theta].
\end{displaymath}

\noindent
The corresponding FC series is then

\begin{displaymath}
  \bar{S}_{\rm c}
  =
  -\,
  \frac{8}{\pi^{2}}
  \sum_{j=0}^{\infty}
  \frac{1}{(2j+1)^{2}}\,
  \sin[(2j+1)\theta],
\end{displaymath}

\noindent
the complex $S_{v}$ series is given by

\begin{displaymath}
  S_{v}
  =
  -\,
  \frac{8}{\pi^{2}}
  \sum_{j=0}^{\infty}
  \frac{1}{(2j+1)^{2}}\,
  v^{2j+1},
\end{displaymath}

\noindent
and the complex power series $S_{z}$ is given by

\begin{displaymath}
  S_{z}
  =
  -\,
  \frac{8}{\pi^{2}}
  \sum_{j=0}^{\infty}
  \frac{1}{(2j+1)^{2}}\,
  z^{2j+1}.
\end{displaymath}

\noindent
The ratio test tells us that the disk of convergence of $S_{z}$ is the
unit disk. If we consider the inner analytic function $w(z)$ within this
disk we observe that $w(0)=0$, as expected. We have for the function and
its derivative

\noindent
\begin{eqnarray*}
  w(z)
  & = &
  -\,
  \frac{8}{\pi^{2}}
  \sum_{j=0}^{\infty}
  \frac{1}{(2j+1)^{2}}\,
  z^{2j+1}
  \;\;\;\Rightarrow
  \\
  \frac{dw(z)}{dz}
  & = &
  -\,
  \frac{8}{\pi^{2}}
  \sum_{j=0}^{\infty}
  \frac{1}{2j+1}\,
  z^{2j}.
\end{eqnarray*}

\noindent
Observe that we have for the derivative the particular value

\begin{displaymath}
  \frac{dw}{dz}(0)
  =
  -\,
  \frac{8}{\pi^{2}}.
\end{displaymath}

\noindent
We may now multiply by $z$ and differentiate again, to obtain

\noindent
\begin{eqnarray*}
  \frac{d}{dz}
  \left[
    z\,
    \frac{dw(z)}{dz}
  \right]
  & = &
  -\,
  \frac{8}{\pi^{2}}\,
  \frac{d}{dz}
  \sum_{j=0}^{\infty}
  \frac{1}{2j+1}\,
  z^{2j+1}
  \\
  & = &
  -\,
  \frac{8}{\pi^{2}}
  \sum_{j=0}^{\infty}
  z^{2j}.
\end{eqnarray*}

\noindent
This is a geometrical series and therefore we may write this expression in
closed form,

\noindent
\begin{eqnarray*}
  \frac{d}{dz}
  \left[
    z\,
    \frac{dw(z)}{dz}
  \right]
  & = &
  -\,
  \frac{8}{\pi^{2}}\,
  \frac{1}{1-z^{2}}
  \\
  & = &
  -\,
  \frac{8}{\pi^{2}}\,
  \left(
    \frac{1}{1+z}
    +
    \frac{1}{1-z}
  \right).
\end{eqnarray*}

\noindent
It is now very simple to integrate in order to obtain the derivative of
$w(z)$ in closed form, remembering that we should have for it the value
$-8/\pi^{2}$ at $z=0$,

\noindent
\begin{eqnarray*}
  z\,
  \frac{dw(z)}{dz}
  & = &
  -\,
  \frac{4}{\pi^{2}}\,
  \left[
    \ln(1+z)
    -
    \ln(1-z)
  \right]
  \\
  & = &
  -\,
  \frac{4}{\pi^{2}}\,
  \ln\!\left(\frac{1+z}{1-z}\right).
\end{eqnarray*}

\noindent
This is the logarithmic derivative of $w(z)$, which is another inner
analytic function, and it should be noted that it is proportional to the
analytic function $w_{\rm sq}(z)$ for the case of the standard square
wave,

\begin{displaymath}
  z\,
  \frac{dw(z)}{dz}
  =
  -\,
  \frac{2}{\pi}\,
  w_{\rm sq}(z).
\end{displaymath}

\noindent
We get for the derivative of $w(z)$

\begin{displaymath}
  \frac{dw(z)}{dz}
  =
  -\,
  \frac{4}{\pi^{2}}\,
  \frac{1}{z}\,
  \ln\!\left(\frac{1+z}{1-z}\right).
\end{displaymath}

\noindent
This function has two borderline hard singularities at the points where
the triangular wave is not differentiable. Presumably $w(z)$ has two
borderline soft singularities at these points. It seems that the second
integration cannot be done explicitly because the indefinite integral of
the function above cannot be expressed as a finite combination of
elementary functions. We are thus unable to write $w(z)$ in closed form.
However, we may still obtain partial confirmation of our results by using
the closed form for the derivative of $w(z)$. If we differentiate the
Fourier series of the triangular wave with respect to $\theta$ we get

\begin{displaymath}
  \frac{dS_{\rm c}}{d\theta}
  =
  \frac{8}{\pi^{2}}
  \sum_{j=0}^{\infty}
  \frac{1}{2j+1}\,
  \sin[(2j+1)\theta].
\end{displaymath}

\noindent
This is equal to $2/\pi$ times the Fourier series for the standard square
wave, given in Subsection~\ref{APPsquarewave}. On the other hand, if we
use $z=\rho\exp(\ii\theta)$ and consider the complex derivative of $w(z)$
taken in the direction of $\theta$, we have

\begin{displaymath}
  \frac{dw(z)}{d\theta}
  =
  \ii
  z\,
  \frac{dw(z)}{dz}.
\end{displaymath}

\noindent
We have therefore

\noindent
\begin{eqnarray*}
  \frac{dw(z)}{d\theta}
  & = &
  -\ii\,
  \frac{4}{\pi^{2}}\,
  \ln\!\left(\frac{1+z}{1-z}\right)
  \\
  & = &
  -\ii\,
  \frac{2}{\pi}\,
  w_{\rm sq}(z).
\end{eqnarray*}

\noindent
The factor of $-\ii$ simply implements the necessary exchanges of real and
imaginary parts, to account for the exchange of sines and cosines in the
process of differentiation, and the factor of $2/\pi$ corrects the
normalization. We see therefore that at least the derivative of the
Fourier series at the unit circle is represented correctly by the inner
analytic function $w(z)$.

\subsection{A Singular Cosine Series}

Consider the Fourier series of the Dirac delta ``function'' centered at
$\theta=\theta_{1}$, which we denote by $\delta(\theta-\theta_{1})$. We
may easily calculate its Fourier coefficients using the rules of
manipulation of $\delta(\theta-\theta_{1})$, thus obtaining
$\alpha_{k}=\cos(k\theta_{1})/\pi$ and $\beta_{k}=\sin(k\theta_{1})/\pi$
for all $k$. The series is therefore the complete Fourier series given by

\noindent
\begin{eqnarray*}
  S_{\rm c}
  & = &
  \frac{1}{2\pi}
  +
  \frac{1}{\pi}
  \sum_{k=1}^{\infty}
  \left[
    \cos(k\theta_{1})
    \cos(k\theta)
    +
    \sin(k\theta_{1})
    \sin(k\theta)
  \right]
  \\
  & = &
  \frac{1}{2\pi}
  +
  \frac{1}{\pi}
  \sum_{k=1}^{\infty}
  \cos(k\Delta\theta),
\end{eqnarray*}

\noindent
where $\Delta\theta=\theta-\theta_{1}$. Apart from the constant term this
is in fact a cosine series on this new variable. Clearly, this series
diverges at all points in the interval $[-\pi,\pi]$. Undaunted by this, we
construct the FC series, with respect to the new variable $\Delta\theta$,

\begin{displaymath}
  \bar{S}_{\rm c}
  =
  \frac{1}{\pi}
  \sum_{k=1}^{\infty}
  \sin(k\Delta\theta),
\end{displaymath}

\noindent
a series that is also divergent, this time almost everywhere. If we define
$v=\exp(\ii\theta)$ and $v_{1}=\exp(\ii\theta_{1})$ the corresponding
complex series $S_{v}$ is then given by

\begin{displaymath}
  S_{v}
  =
  \frac{1}{2\pi}
  +
  \frac{1}{\pi}
  \sum_{k=1}^{\infty}
  \left(\frac{v}{v_{1}}\right)^{k},
\end{displaymath}

\noindent
where we included the $k=0$ term, and the corresponding complex power
series $S_{z}$ is given by

\begin{displaymath}
  S_{z}
  =
  \frac{1}{2\pi}
  +
  \frac{1}{\pi}
  \sum_{k=1}^{\infty}
  \left(\frac{z}{z_{1}}\right)^{k},
\end{displaymath}

\noindent
where $z=\rho v$ and $z_{1}=v_{1}$ is a point on the unit circle. This is
a geometrical series that converges everywhere strictly inside the unit
disk, where it therefore converges to the analytic function

\begin{displaymath}
  w_{\delta}(z)
  =
  \frac{1}{2\pi}
  -
  \frac{1}{\pi}\,
  \frac{z}{z-z_{1}}.
\end{displaymath}

\noindent
This function has a simple pole at the singular point of the Dirac delta
``function''. This is a hard singularity of degree $1$. As was discussed
in Section~\ref{SECcomprel}, strictly speaking this is not an inner
analytic function by the definition which we adopted here, for two
reasons. First, $w_{\delta}(0)$ is not zero, which is easily fixed by just
taking off the constant term. Second, it does not reduce to a real
function over the real axis. However, as was shown in
Section~\ref{SECcomprel}, it does reduce to a real function over the
straight line $z=\chi z_{1}$, with real $\chi$, since in this case we have

\begin{displaymath}
  w(z)
  =
  -\,
  \frac{1}{\pi}\,
  \frac{\chi}{\chi-1},
\end{displaymath}

\noindent
which is a real function of the real variable $\chi$. Since we have
$z_{1}=\exp(\ii\theta_{1})$, we see that this is an inner analytic
function rotated by an angle $\theta_{1}$ around the origin. We see
therefore that our definition of inner analytic function can be easily
generalized in this way, to make reference to any straight line going
through the origin. We also see that there is nothing to prevent us from
treating this function just like any other inner analytic function.

We now examine the real part of the function $w_{\delta}(z)$, which was
derived in Section~\ref{SECcomprel}, and which is related to the series
$S_{\rm c}$,

\begin{displaymath}
  \Re[w_{\delta}(z)]
  =
  \frac{1}{2\pi}
  -
  \frac{1}{\pi}\,
  \frac
  {\rho\left[\rho-\cos(\Delta\theta)\right]}
  {\left(\rho^{2}+1\right)-2\rho\cos(\Delta\theta)}.
\end{displaymath}

\noindent
The $\rho\to 1$ limit of $\Re[w_{\delta}(z)]$ gives us back the constant
term and the original cosine series $S_{\rm c}$, and in fact attributes
well-defined finite values to it almost everywhere, even though the
original series is divergent. As was shown directly in
Section~\ref{SECcomprel}, in the $\rho\to 1$ limit this real part has all
the required properties of the Dirac delta ``function'', including the
fact that it assumes the value zero almost everywhere. Note that although
the singularity of $w_{\delta}(z)$ at $z=z_{1}$ is a simple pole and thus
not an integrable one, the real part of $w_{\delta}(z)$ is integrable if
one crosses the singularity in a specific direction, which in this case is
the direction of the integration along the unit circle.

It is interesting that we may also write a closed expression for the real
function which is the FC function of the Dirac delta ``function''. We just
consider the imaginary part of $w_{\delta}(z)$, which is related to the
series $\bar{S}_{\rm c}$,

\begin{displaymath}
  \Im[w_{\delta}(z)]
  =
  \frac{1}{\pi}\,
  \frac
  {\rho\sin(\Delta\theta)}
  {\left(\rho^{2}+1\right)-2\rho\cos(\Delta\theta)}.
\end{displaymath}

\noindent
We now take the $\rho\to 1$ limit, assuming that $\Delta\theta\neq 0$, and
since this is an actual function we may write

\noindent
\begin{eqnarray*}
  \bar{f}_{\rm c}(\Delta\theta)
  & = &
  \lim_{\rho\to 1}
  \Im[w_{\delta}(z)]
  \\
  & = &
  \frac{1}{2\pi}\,
  \frac{\sin(\Delta\theta)}{1-\cos(\Delta\theta)}
  \\
  & = &
  \frac{1}{2\pi}\,
  \frac{1+\cos(\Delta\theta)}{\sin(\Delta\theta)}.
\end{eqnarray*}

\noindent
This is just a regular function except for a simple pole at
$\Delta\theta=0$, and in fact quite a simple rational function involving
trigonometric functions.

Note that we get a free set of integration formulas out of this effort.
Since the coefficients $a_{k}$ can be written as integrals involving
$\bar{f}_{\rm c}(\Delta\theta)$, we have at once that for all $k>0$

\begin{displaymath}
  \int_{-\pi}^{\pi}d\theta\,
  \frac{1+\cos(\Delta\theta)}{\sin(\Delta\theta)}\,
  \sin(k\Delta\theta)
  =
  2\pi.
\end{displaymath}

\noindent
We may easily construct from this single Dirac delta ``function'' a pair
of such ``functions'' which relates to the derivative of the standard
square wave. All we have to do is to use the values $\theta_{1}=0$ and
$\theta_{1}=\pi$ and add to the delta ``function'' at zero given by
$2\delta(\theta-0)$ the delta ``function'' at $\pi$ given by
$-2\delta(\theta-\pi)$. If we call the corresponding analytic functions
$w_{0}(z)$ and $w_{\pi}(z)$, we have for the function $w_{2\delta}(z)$
corresponding to the pair

\noindent
\begin{eqnarray*}
  w_{0}(z)
  & = &
  \frac{1}{\pi}
  +
  \frac{2}{\pi}\,
  \frac{z}{1-z}
  \;\;\;\Rightarrow
  \\
  w_{\pi}(z)
  & = &
  -\,
  \frac{1}{\pi}
  +
  \frac{2}{\pi}\,
  \frac{z}{1+z}
  \;\;\;\Rightarrow
  \\
  w_{2\delta}(z)
  & = &
  \frac{2}{\pi}\,
  \frac{z}{1-z}
  +
  \frac{2}{\pi}\,
  \frac{z}{1+z}
  \\
  & = &
  \frac{4}{\pi}\,
  \frac{z}{1-z^{2}}.
\end{eqnarray*}

\noindent
This is an inner analytic function with two simple poles, one at $z=1$ and
one at $z=-1$. If, on the other hand, we consider the logarithmic
derivative of the analytic function $w_{\rm sq}(z)$ of the standard square
wave, we get

\noindent
\begin{eqnarray*}
  z\,
  \frac{dw_{\rm sq}(z)}{dz}
  & = &
  \frac{2}{\pi}\,
  z\,
  \frac{d}{dz}
  \ln\!\left(\frac{1+z}{1-z}\right)
  \\
  & = &
  \frac{2}{\pi}\,
  z\,
  \left(
    \frac{1}{1+z}
    +
    \frac{1}{1-z}
  \right)
  \\
  & = &
  \frac{4}{\pi}\,
  \frac{z}{1-z^{2}}.
\end{eqnarray*}

\noindent
We see therefore that we do indeed have

\begin{displaymath}
  w_{2\delta}(z)
  =
  z\,
  \frac{dw_{\rm sq}(z)}{dz}.
\end{displaymath}

\noindent
In terms of derivatives with respect to $\theta$ this is written as

\begin{displaymath}
  w_{2\delta}(z)
  =
  -\ii\,
  \frac{dw_{\rm sq}(z)}{d\theta},
\end{displaymath}

\noindent
where once again the factor of $-\ii$ has the role of exchanging real and
imaginary parts, and thus sines and cosines.

Observe that this is an example in which, although the function
$w_{\delta}(z)$ is analytic on the open unit disk, the series $S_{z}$ is
not convergent on any points of the unit circle, and hence the two FC
Fourier series do not converge together on any points of the unit circle.
It is therefore outside the hypotheses we used in most of this paper, but
it is still full of meaning. This indicates that there is more in this
structure than has been analyzed so far. If fact, it is quite possible
that the whole structure of distributions such as the one related to the
Dirac delta ``function'' is lurking on the rim of the maximum disk of
convergence of the series $S_{z}$ of the inner analytic functions
$w_{\delta}(z)$ corresponding to divergent trigonometric series, that is,
on the rim of the unit disk.

\end{document}